\documentclass[a4paper]{article}
\usepackage[mathscr]{eucal}
\usepackage{amsmath,amsfonts,amssymb,amsthm,amscd}
\usepackage{graphicx}
\usepackage[square]{natbib}
\usepackage{multirow}
\newcommand\field{\mathbb}

\newcommand\R{\field{R}}

\newcommand\C{\field{C}}
\newcommand\Z{\field{Z}}
\newcommand\M{\field{M}}

\newcommand\Q{\field{Q}}

\newcommand\tr{\operatorname{Tr}}

\newcommand\grad{\operatorname{grad}}

\newcommand\card{\operatorname{card}}

\newcommand\id{\operatorname{\mathrm{Id}}}

\newcommand\rmi{\mathrm{i}\mspace{1mu}}

\newcommand\pder[2]{\dfrac{\partial #1 }{\partial #2}}
\newcommand\bfi[1]{{\bfseries\itshape{#1}}}

\newcommand\mtext[1]{\quad\text{#1}\quad}

\newcommand\defset[2]{\left\{{#1}\;\vert \;\; {#2} \,\right\}}

\usepackage[square]{natbib}
\theoremstyle{plain}
\newtheorem{theorem}{Theorem}
\newtheorem{lemma}[theorem]{Lemma}
\newtheorem{proposition}[theorem]{Proposition}
\newtheorem{corollary}[theorem]{Corollary}

\newtheorem{definition}[theorem]{Definition}
\newtheoremstyle{note}{\topsep}{\topsep}{\slshape}{}{\scshape}{}{ }{}
\theoremstyle{note}

\newtheorem{remark}[theorem]{Remark}
\numberwithin{equation}{section}
\numberwithin{theorem}{section}

\begin{document}
\thispagestyle{empty}
\vspace*{1em}
\begin{center}
\LARGE\textbf{Jordan obstruction to the integrability of Hamiltonian systems with homogeneous
potentials}
\end{center}
\vspace*{0.5em}
\begin{center}
\large  Guillaume Duval$^1$ and Andrzej J.~Maciejewski$^2$
\end{center}
\vspace{2em}
\hspace*{2em}\begin{minipage}{0.8\textwidth}
\small
$^1$1 Chemin du Chateau , 76 430 Les Trois Pierres, France. \\
(e-mail: dduuvvaall@wanadoo.fr)\\[0.5em]
$^2$Institute of Astronomy,
  University of Zielona G\'ora,
  Podg\'orna 50, \\\quad PL-65--246 Zielona G\'ora, Poland
  (e-mail: maciejka@astro.ia.uz.zgora.pl)
\end{minipage}

\begin{abstract}
In this paper, we consider the natural complex Hamiltonian systems with homogeneous potential $V(q)$, $q\in\C^n$,
of degree $k\in\Z^\star$. The known results of Morales and Ramis give necessary conditions for the complete integrability of such systems. These conditions are expressed in terms of the eigenvalues of the Hessian matrix $V''(c)$ calculated at a non-zero point $c\in\C^n$, such that $V'(c)=c$. The main aim of this paper is to show that there are other obstructions for the integrability which appear if the  matrix $V''(c)$ is not diagonalizable. We prove, among other things, that if  $V''(c)$ contains a Jordan block of size greater  than  two, then the system is not integrable in the Liouville sense. The main ingredient in the proof of this result consists in translating  some ideas of Kronecker about Abelian extensions of number fields into the framework of  differential Galois theory.
\end{abstract}
\section{Introduction}

\subsection{Morales and Ramis results}

The Galois obstruction to the integrability of Hamiltonian systems is formulated in the
following theorem obtained by Morales and Ramis~\cite{Morales:01::b1}.

\begin{theorem}[Morales-Ramis]
  If an Hamiltonian system is completely integrable with first
  integrals  meromorphic in a connected neighbourhood of a phase curve $\gamma$, then the identity
  component of the differential Galois group of the variational equation along
  $\gamma$ is virtually Abelian.
\end{theorem}
In \cite{moralesramis01}, Morales and Ramis applied this theorem to find
obstructions to the complete integrability of Hamiltonian systems with
homogeneous potentials. They considered natural systems  with Hamiltonian given by
\begin{equation}
  H (q, p) = \frac{1}{2}  \sum_{i = 1}^n p_i^2 + V (q_1, \ldots, q_n),
  \label{eqH}
\end{equation}
where $q=  (q_1, \ldots, q_n), p=(p_1,\ldots, p_n)\in\C^n$, are the canonical coordinates and momenta, respectively, and  $V (q)$ is a homogeneous potential of degree $k \in \Z^\star:= \Z
\setminus \{0\}$.
To that purpose, following Yoshida~\cite{yoshida87}, they studied \bfi{the
variational equations} (in short: VE) associated to \bfi{a
proper Darboux point} of $V$, (in short: PDP) which is a non-zero vector $c\in\C^n$ such that
\begin{equation}
 \grad V(c)=: V' (c) = c. \label{eqdarboux}
\end{equation}
If such a Darboux point exists, then  the Hamiltonian
system admits  a particular solution associated  with this point, namely,  the rectilinear trajectory
\begin{equation*}
  t \mapsto \gamma (t) =(q(t), p(t)) :=( \varphi(t) c, \dot\varphi(t) c)
\in\C^{2 n},
\end{equation*}
where $t \mapsto \varphi (t)$ is a complex scalar function satisfying the
hyper-elliptic differential equation
\begin{equation}
  \dot\varphi (t)^2 = \frac{2}{k} (1 - \varphi^k (t)) \ \Longrightarrow\  \ddot\varphi (t) =
  - \varphi^{k - 1} (t) . \label{eqelliptic}
\end{equation}
The VE along the curve $t \mapsto \gamma (t)$ is given  by
\begin{equation}
  \frac{d^2 \eta}{dt^2} = - \varphi^{k - 2} (t) V'' (c) \eta, \quad \eta\in \C^n. \label{eqVEt}
\end{equation}
The Hessian matrix $V'' (c)$ is a $n \times n$ complex, symmetric
scalar matrix. Assume that it is diagonalizable with eigenvalues
$(\lambda_1, \ldots, \lambda_n) \in\C^n$ (we called them \bfi{the Yoshida coefficients}).  Then up to a linear change
of unknowns, the system~\eqref{eqVEt} splits into a
direct sum of equations
\begin{equation}
  \frac{d^2 \eta_i}{dt^2} = - \lambda_i \varphi^{k - 2} (t) \eta_i, \quad i=1,\ldots,n.
  \label{eqVElambda}
\end{equation}
\begin{table}[th]
\caption{The Morales-Ramis table. }
\begin{center}
  \begin{tabular}{|c|c|c|c|}
\hline
\vbox to 1em {}    $\mathcal{G}(k, \lambda)^{\circ}$ & $k$ & $\lambda$ & Row number\\
    \hline
    & $k = \pm 2$ & $\lambda$ is an arbitrary complex number & 1\\
    \hline
\vbox to 1.6em {}    $G_{\mathrm{a}}$ & $|k| \geq 3$ & $\lambda (k , p) = p + \dfrac{k}{2} p (p - 1)$
    & 2\\[0.5em]
    \hline
\vbox to 1.6em {}   $G_{\mathrm{a}}$  & 1 & $p + \dfrac{1}{2} p (p - 1)$, $p \neq - 1 , 0$ & 3\\[0.5em]
    \hline
\vbox to 1.6em {}   $G_{\mathrm{a}}$   & -1 & $p - \dfrac{1}{2} p (p - 1)$, $p \neq 1 , 2$ & 4\\[0.5em]
    \hline
    $\{\id\}$   & 1 & $0$ & 5\\
    \hline
  $\{\id\}$   & -1 & $1$ & 6\\
    \hline
\vbox to 1.6em {}  $\{\id\}$     & $|k| \geq 3$ & $\dfrac{1}{2}  \left( \dfrac{k - 1}{k} + p (p + 1) k\right)
   $ & 7\\[0.5em]
    \hline
\vbox to 1.6em {}  \multirow{2}{*}{ $\{\id\}$  }  & \multirow{2}{*}{3} & $\dfrac{- 1}{24} + \dfrac{1}{6} (1 + 3 p)^2$, $\dfrac{- 1}{24} +
    \dfrac{3}{32} (1 + 4 p)^2$ & 8,9\\[0.5em]
&    & $\dfrac{- 1}{24} + \dfrac{3}{50} (1 + 5 p)^2$, $\dfrac{- 1}{24} +
    \dfrac{3}{50} (2 + 5 p)^2$ & 10,11\\[0.5em]
   \hline
\vbox to 1.6em {}    \multirow{2}{*}{ $\{\id\}$  }   &  \multirow{2}{*}{-3}& $\dfrac{25}{24} - \dfrac{1}{6} (1 + 3 p)^2$, $\dfrac{25}{24} -
    \dfrac{3}{32} (1 + 4 p)^2$ & 12,13\\[0.5em]
\vbox to 1.6em {}     &  & $\dfrac{25}{24} - \dfrac{3}{50} (1 + 5 p)^2$, $\dfrac{25}{24} -
    \dfrac{3}{50} (2 + 5 p)^2$ & 14,15\\[0.5em]
    \hline
\vbox to 1.6em {}   $\{\id\}$    & 4 & $\dfrac{- 1}{8} + \dfrac{2}{9} (1 + 3 p)^2$ & 16\\[0.5em]
    \hline
\vbox to 1.6em {}   $\{\id\}$    & -4 & $\dfrac{9}{8} - \dfrac{2}{9} (1 + 3 p)^2$ & 17\\[0.5em]
    \hline
 \vbox to 1.6em {}  $\{\id\}$   & 5 & $\dfrac{- 9}{40} + \dfrac{5}{18} (1 + 3 p)^2$, $\dfrac{- 9}{40} +
    \dfrac{1}{10} (2 + 5 p)^2$ & 18,19\\[0.5em]
    \hline
  \vbox to 1.6em {}  $\{\id\}$   & -5 & $\dfrac{49}{40} - \dfrac{5}{18} (1 + 3 p)^2$, $\dfrac{49}{40} -
    \dfrac{1}{10} (2 + 5 p)^2$ & 20,21\\[0.5em]
    \hline
  \end{tabular}
\end{center}
\label{tableMR}
\end{table}

Morales and Ramis proved the following
\begin{theorem}[Morales-Ramis]
  \label{MR} Assume that the Hamiltonian system with Hamiltonian \eqref{eqH} and
  $\deg (V) = k \in \Z^\star$ is completely integrable by meromorphic first
  integrals. If  $c = V' (c)$ is a PDP of $V$ and the Hessian
  matrix $V'' (c)$ is diagonalizable with the Yoshida coefficients
  $(\lambda_1, \ldots \lambda_n) \in\C^n$, then each pair $(k ,
  \lambda_i)$ belongs to Table~\ref{tableMR}.
\end{theorem}
The group $\mathcal{G}(k, \lambda)^{\circ}$ appearing in the first column of Table~\ref{tableMR} will be properly
defined in Section~1.3.

\subsection{Jordan obstruction.}

In order to generalise Theorem \ref{MR}, we are going to work without  the  assumption  that the Hessian matrix $V'' (c)$ is semi-simple.  Indeed, since the Hessian matrix $V''
(c)$ is symmetric,  it is diagonalizable if it is real. But,  even
for a real potential coming from physics, PDP may be a complex non real vector.
Therefore, $V'' (c)$ may not be
diagonalizable, see Section 6 for a discussion about this point.

Our main result is the following.
\begin{theorem}
  \label{abclassification} Let $V (q)$ be a homogeneous
  potential of $n$ variables and degree $k \in \Z \setminus \{- 2,
  0, 2\}$, such that $H$ is completely integrable with meromorphic first
  integrals. Then, at any proper Darboux point $c = V' (c) \in\C^n
  \setminus \{0\}$, the Hessian matrix $V'' (c)$ satisfies the
  following conditions:
  \begin{enumerate}
    \item For each eigenvalue $\lambda$ of $V'' (c)$, the pair $(k, \lambda)$
    belongs to Table \ref{tableMR}.
    \item The matrix $V'' (c)$ does not have any Jordan block of size $d \geq 3$.
    \item If $V'' (c)$ has a Jordan block of size $d = 2$  with corresponding eigenvalue $\lambda$, then  the row number of $(k, \lambda)$ in Table \ref{tableMR} is
    greater or equal to five.
  \end{enumerate}
 For  $k = \pm 2$, independently of the value of $V'' (c)$, the connected
  component of the Galois group of the variational equation is Abelian.
\end{theorem}
In the above statement,  by \bfi{a  Jordan block of size $d$ with the  eigenvalue $\lambda$}, we mean that the Jordan form  of  $V''(c)$ contains a block of the form
\begin{equation}
B(\lambda, d):= \begin{bmatrix}
                \lambda & 0 &0& \hdotsfor{2}& 0 \\
                 1& \lambda & 0& \hdotsfor{2}&0\\
                 0& 1&\lambda& \hdotsfor{2}&0 \\
                    \hdotsfor{6}\\
                 0& 0& 0& \ldots &1& \lambda
               \end{bmatrix}
\in\M(d,\C),
\end{equation}
where $\M(d,\C)$ denotes the set of $d\times d$ complex matrices.

\begin{remark}
  Theorem \ref{abclassification} roughly states that Morales-Ramis Theorem
  \ref{MR} is optimal. Indeed, up to some exceptions, if $H$ is completely
  integrable, then $V'' (c)$ must be diagonalizable with specific eigenvalues.

  Our result  is  analogous  to the
  Liapunov-Kowaleskaya Theorem, which states that if a given system of weight-homogeneous differential equations enjoys  the Painleve property, then among other things, the linearization of the system along a certain single-valued particular solution is diagonalizable. For details, see \cite{Kozlov:96::}. Moreover, in the same sense, we find similarities in the classical normal form theory of vector fields, where a complicated dynamics appears in a neighbourhood of the equilibrium if the linearization of the vector field is not semi-simple.

As far as we know, except for one example given in Chapter~7 in \cite{Morales:99::c}, there are no explicit links between the Galois approach to the integrability and the dynamics. Nevertheless, the above analogies were our strong motivations for that study.

  The proof of Theorem \ref{abclassification} is of another nature. It comes
  from arithmetic ideas belonging to Kronecker. He
  observed that in  Number Theory,  Abelian extensions of number fields can be characterised by
  simple arithmetic relations. We  translate this very nice observation  into the framework of the Differential Galois Theory.

\end{remark}

\subsection{VE, Yoshida transformations and Jordan blocks}

The VE \eqref{eqVEt} is a system of differential equations with respect to the time
variable $t$. First,  we perform   the so-called \bfi{Yoshida transformation}, in
order to express the VE in terms of a new variable $z$. The great advantage of
this transformation is that it  converts our original system into a new one where the
classical hypergeometric equation appears naturally. Next, we give the
canonical formulae for the subsystems of VE  associated to  Jordan
blocks.

The Yoshida transformation  is  a change of   independent variable in  equation \eqref{eqVEt}
 given by
\begin{equation}
t \longmapsto  z = \varphi^k (t) \label{eqyoshida} .
\end{equation}
Thanks to \eqref{eqelliptic} and the chain rule we have
\begin{gather*}
 \frac{d^2 \eta}{dt^2} =\left (
\frac{dz}{dt}\right)^2  \frac{d^2 \eta}{dz^2} + \frac{d^2 z}{dt^2}  \frac{d
\eta}{dz}, \\
\left(\frac{dz}{dt}\right)^2 =2kz(1-z)\varphi^{k-2}(t), \quad  \frac{d^2 z}{dt^2}=[ (2-3k)z+2(k-1)]\varphi^{k-2}(t).
\end{gather*}
Then, after some simplifications, \eqref{eqVEt}, becomes
\begin{equation}
 \label{eqVEz}
\frac{d^2 \eta}{dz^2} + p (z) \frac{d \eta}{dz}  =  s (z) V'' (c) \eta,
\end{equation}
where
\begin{equation*}
 p (z) = \frac{2 (k - 1) (z - 1) + kz}{2 kz (z - 1)} \mtext{ and } s (z) =
   \frac{1}{2 kz (z - 1)} .
\end{equation*}
%
Next, after the classical Tchirnhauss change of dependent variable,
\begin{equation}
\label{eqtchirn}
  \eta = f (z) \zeta,  \qquad f (z) = \exp\Big( -\frac{1}{2} \int p (z) dz\Big) =
  z^{-( k - 1)/2 k} (z - 1)^{-1/4},
\end{equation}
 equation \eqref{eqVEz} has the  reduced
form
\begin{equation}
\label{eqVEredz}
  \frac{d^2 \zeta}{dz^2} = [ r_0 (z) \id + s (z) V'' (c)] \zeta,
\end{equation}
where
\[ r_0 (z) = \frac{\rho^2 - 1}{4 z^2} + \frac{\sigma^2 - 1}{4 (z - 1)^2} -
   \frac{1}{4}\left (1 - \rho^2 - \sigma^2 + \tau_0^2\right)\left( \frac{1}{z} + \frac{1}{1 -
   z}\right),
\]
and
\begin{equation*}
 \rho = \frac{1}{k}, \qquad  \sigma = \frac{1}{2}, \qquad \tau_0 = \frac{k - 2}{2 k}.
\end{equation*}

Assume that  $V'' (c)$ contains a Jordan block $B(\lambda,d)$ with $d=3$, for example.
Then, the subsystem of \eqref{eqVEredz} corresponding to this block can be written as
\begin{equation*}
\frac{d \phantom{z}}{d z^2} \begin{bmatrix}
 x  \\ y\\  u
\end{bmatrix}=
 \begin{bmatrix}
 x''  \\ y''\\  u''
\end{bmatrix}=
\left(
\begin{bmatrix}
  r_0 (z) & 0 & 0\\
    0 & r_0 (z) & 0\\
    0 & 0 & r_0 (z)
\end{bmatrix}
+ s(z)
\begin{bmatrix}
   \lambda & 0 & 0\\
    1 & \lambda & 0\\
    0 & 1 & \lambda
\end{bmatrix}
\right)
\begin{bmatrix}
 x\\ y\\ u
\end{bmatrix}.
\end{equation*}
We  rewrite it in the following form
\begin{equation}
 \label{eqblock}
  \begin{bmatrix}
 x''  \\ y''\\  u''
\end{bmatrix}=
\begin{bmatrix}
  r (z) & 0 & 0\\
    s (z) & r_{} (z) & 0\\
    0 & s (z) & r (z)
\end{bmatrix}
\begin{bmatrix}
 x\\ y\\ u
\end{bmatrix}.
\end{equation}
where $r (z) = r_{\lambda} (z)$ is given by
\begin{equation}
\label{eq:rl}
 r (z) = r_0 (z) + \lambda s (z) = \frac{\rho^2 - 1}{4 z^2} + \frac{\sigma^2
   - 1}{4 (z - 1)^2} - \frac{1}{4} \left(1 - \rho^2 - \sigma^2 + \tau^2\right)\left (
   \frac{1}{z} + \frac{1}{1 - z}\right),
\end{equation}
with
\begin{equation}
 \label{eqdeltaexp}
  \rho = \frac{1}{k}, \qquad \sigma = \frac{1}{2},  \qquad \tau = \frac{\sqrt{(k - 2)^2 + 8 k
  \lambda}}{2 k} .
\end{equation}
The above three numbers are exactly  the respective exponents differences at $z
= 0$,  $z = 1$ and  $z = \infty$ of the reduced hypergeometric
equation $L_2 = x'' - r (z) x = 0$. Thus,  the solutions of $L_2 = 0$ belong to the
Riemann scheme
\begin{equation}
  P\left\{ \begin{matrix}
    0 & 1 & \infty\\
    \dfrac{1}{2} - \dfrac{1}{2 k} & \dfrac{1}{4} & \dfrac{- 1 - \tau}{2}\\[0.5em]
    \dfrac{1}{2} + \dfrac{1}{2 k} & \dfrac{3}{4} & \dfrac{- 1 + \tau}{2}
  \end{matrix}\; z \right\} \label{eqRS}.
\end{equation}
 The group $\mathcal{G}(k, \lambda)^\circ$ appearing in the first column of Table~\ref{tableMR}, is the identity component of the differential Galois group of the equation $L_2
= x'' - r (z) x = 0$,  with respect  to the ground field $\C(z)$.

\subsection{Generalities and Galois groups of the distinct VE}

In this subsection, we summarise some results about differential Galois groups
and classical Differential Algebra which  we frequently use, see \cite{kolchin68,svdp}. Next, we compare
the differential Galois groups of  different forms of  the VE
introduced in the previous subsection.

In what follows, $(K,\partial)$ denotes an ordinary differential field with the algebraically closed subfield of constants $\C$. We use the standard notation, e.g.,  $x'=\partial x$,  $x''=\partial^2 x$, etc., for an element $x\in K$.

\begin{itemize}
  \item If a linear system $Y' = AY$, where $A\in\M(n_1+n_2,K)$,  splits into a direct sum
\begin{equation*}
Y'=
 \begin{bmatrix}
      Y_1'\\
       Y_2'
\end{bmatrix}
=
 \begin{bmatrix}
     A_1 & 0\\
       0 & A_2
\end{bmatrix}
 \begin{bmatrix}
      Y_1\\
       Y_2
\end{bmatrix}
\mtext{where} A_i\in\M(n_i,K) \mtext{for}i=1,2,
\end{equation*}
  then, with obvious notations,  the identity component $G^{\circ}$ of its differential Galois group $G$ is a subgroup of the direct product
  $G_1^{\circ} \times G_2^{\circ}$. Moreover, the two  projection maps $\pi_i :
  G^{\circ} \rightarrow G_i^{\circ}$, with $i=1,2$, are  surjective. Therefore, $G^{\circ}$ is
  Abelian iff $G_1^{\circ}$ and $G_2^{\circ}$ are Abelian.
  \item If the system $Y_1' = A_1 Y_1$ is a subsystem of
\begin{equation*}
 \begin{bmatrix}
      Y_1'\\
       Y_2'
\end{bmatrix}
=
 \begin{bmatrix}
     A_1 & 0\\
       B  & A_2
\end{bmatrix}
 \begin{bmatrix}
      Y_1\\
       Y_2
\end{bmatrix},
\end{equation*}
then the reduction morphism $G^{\circ} \rightarrow
  G_1^{\circ}$ is surjective. Therefore, if $G^{\circ}$ is  Abelian, then
  $G_1^{\circ}$ is also Abelian.
\end{itemize}
\begin{lemma}
  \label{ostrowski}Let $E / K$ be an ordinary differential field
  extension with the same subfield of constants $\C$.
  \begin{enumerate}
    \item Let $f_1, \ldots, f_p\in E$,  and $f_i' \in K$, for $i=1,\ldots, p$.
     Then the family $\{f_1, \ldots, f_p \}$ is algebraically dependent
    over $K$ if and only if there exists a non trivial linear relation
    \[ c_1 f_1 + \cdots + c_p f_p \in K \text{ with } (c_1, \ldots, c_p) \in
       \C^p \backslash \{0\}. \]
    \item Let $T (E / K)$ be the set of elements $f$ of $E$ such that there exists a non-zero linear differential equation $L\in K[\partial]$ such that $L(f)=0$. Then $T (E /
    K)$ is a $K$-algebra containing the algebraic closure of $K$ in $E$.  If
    $E / K$ is a Picard-Vessiot extension, then $T (E / K)$ is the Picard-Vessiot
    ring of $E / K$,  and $T (E / K) = K [Z_{i, j}] [ {W}^{-1}]$, where $W=\det((Z_{i, j}))$, and  $(Z_{i, j})_{1 \leq i, j \leq n}$ is an arbitrary fundamental matrix
    defining the Picard-Vessiot extension.
    \item Let $y' = Ay$ be a differential system with $A \in \M(n,K)$, and  $K'
    / K$ be a finite degree  extension of $K$. Denote by $G$
    (resp. by $G'$) the respective Galois groups of $y' = Ay$ when this system is
    considered over $K$, (resp. over $K'$). Then $G'$ is naturally a subgroup
    of $G$,  and  $(G')^{\circ} = G^{\circ}$.
  \end{enumerate}
\end{lemma}
\begin{proof}
  (1)  is  the classical Ostrowski-Kolchin theorem about the algebraic independence of
  integrals. Its proof may be found in {\cite{kolchin68}}.

\noindent
(2) follows directly from
   Exercises~1.24 on  p. 17 and  Corollary~1.38 on  p. 30 in~\cite{svdp}.

\noindent
  (3)  Let $F' / K'$ be a Picard-Vessiot extension of $y' = Ay$ over $K'$.
  Then $F' = K' (Y)$, where $Y$ is a fundamental matrix of solutions of the
  system. Set $F = K (Y)$. Then  $F / K$ is a Picard-Vessiot extension of $y' =
  Ay$ over $K$, for which
  \[ F' = K' F \mtext{ and } F \subset F' . \]
  Since the group $G'$ fixes $K'$ pointwise and leaves $F$ globally invariant, i.e., $G'\cdot F =F$, it may be
  considered as a subgroup of $G$. Therefore, we also have the inclusion of
  connected components
  \[ (G')^{\circ} \subset G^{\circ} . \]
  From  Corollary 1.30 on p.~23 in \cite{svdp} and its proof, we have
  \[ \dim G^{\circ} = \dim G = \mathrm{tr} . \deg (F / K) . \]
  But $\mathrm{tr} . \deg (F / K) = \mathrm{tr} . \deg (K' F / K') = \mathrm{tr} .
  \deg (F' / K') = \dim (G')^{\circ}$. So  the two irreducible varieties
  ($G')^{\circ}$ and $G^{\circ}$ have the same dimension, hence they are
  equal.
\end{proof}

 Consider the four             variational equations   derived in the previous subsection,  namely  equations
\eqref{eqVEt},  \eqref{eqVEz},  \eqref{eqVEredz}  and  \eqref{eqblock}. The system
\eqref{eqVEt} is defined over the ground fields $\C(\varphi (t),
\dot\varphi (t))$ and the other three are defined over $\C(z)$. Let
 $G (\mathrm{VE}_t)$,   $G (\mathrm{VE}_z)$,
and $G_{\mathrm{block}}$  be the differential  Galois groups of  equations
\eqref{eqVEt},  \eqref{eqVEz},  and  \eqref{eqblock}, respectively. Then we have   the following.
\begin{proposition}
  \label{VEtVEz}With the notations above, we have:
  \begin{enumerate}
    \item The Galois groups of systems \eqref{eqVEz} and \eqref{eqVEredz} have common connected
    component $G (\mathrm{VE}_z)^{\circ}$.

    \item The two connected components $G (\mathrm{VE}_t)^{\circ}$ and $G
    (\mathrm{VE}_z)^{\circ}$ are isomorphic.

    \item The connected component $G_{\mathrm{block}}^{\circ}$ is a quotient of
    $G (\mathrm{VE}_z)^{\circ}$.
  \end{enumerate}
\end{proposition}
\begin{proof}
(1)~Set  $K =\C(z)$ and $K' =\C(z) [f (z)]$, where $f
  (z) $ is given by
  \eqref{eqtchirn}. Then $K' / K$ is a finite extension. Denote by
  $Z$ a fundamental matrix of solutions of \eqref{eqVEredz}. Then $f (z) Z
  $ is a fundamental matrix of solutions of \eqref{eqVEz}. Therefore,
   \eqref{eqVEz}  and \eqref{eqVEredz} share the same
  Picard-Vessiot extension over $K'$. So they have the same Galois group $G'$ over $K'$.
  From point 3 of Lemma~\ref{ostrowski},
  \[ (G')^{\circ} = G^{\circ} = G (\mathrm{VE}_z)^{\circ}, \]
  is also the connected component of the Galois group of \eqref{eqVEredz}  when
  it is viewed as a system over $K =\C(z)$.

 (2)~Consider the Yoshida map
\begin{equation*}
 \begin{split}
 \phi :  K=\C(z) & \rightarrow   K'
        =\C(\varphi (t), \dot\varphi (t)),\\
  z&\longmapsto  \varphi^k (t).
\end{split}
\end{equation*}
  This map  is a morphism of fields which is not a differential morphism for  differential fields $(K,\frac{d}{dz})$ and $(K',\frac{d}{dt})$. But,
  since $K' / K$ is finite, the derivation $\frac{d}{dz}  $ of $K$ extends
  uniquely to a derivation of $K'$ which is still denoted by the same symbol.
  Moreover,
  \begin{equation}
    \frac{d}{dt} = \frac{dz}{dt}  \frac{d}{dz} = k \varphi^{k - 1} \dot\varphi
    \frac{d}{dz} \label{chainrule}
  \end{equation}
  Let $F' / K'$ be the Picard-Vessiot extension of \eqref{eqVEt}  over $(K', \frac{d}{dt})$.  Then  $F' / K'$ is a Picard-Vessiot extension
  of \eqref{eqVEz} when considered over $(K', \frac{d}{dz})$. From
  \eqref{chainrule}, an automorphism of $F' / K'$ commutes with $\frac{d}{dt} $
  iff it commutes with $ \frac{d}{dz}$. Therefore, by point 3 of Lemma \ref{ostrowski},  we conclude that
  \[ G (\mathrm{VE}_t)^{\circ} = G (\mathrm{VE}_z)^{\circ}.  \]

  (3)~This point has already been  proved.
\end{proof}

\subsection{The plan of the paper}

As  shown in the above section, there are several  VE,  but  essentially we have
two connected Galois group to deal with: $G (\mathrm{VE}_z)^\circ$ and $G_{\mathrm{block}}^{\circ}$,
the latter being
a quotient of the former.

In Section~2, we study differential equations of the form \eqref{eqblock} for
Jordan blocks of size $d = 2$. We find necessary and sufficient conditions for
the connected component of the Galois group $G_{\mathrm{block}}^{\circ} =
G_2^{\circ}$ to be Abelian, see Theorem~\ref{abcritere}. In this part the reader will find  our interpretation of  Kronecker's ideas  in the framework of the
Differential Galois Theory.

In Section~3, we apply this result to eliminate from Table~\ref{tableMR} all
the cases corresponding to $G_1^{\circ} \simeq G_{\mathrm{a}}$, where $G_{\mathrm{a}}$ denotes the additive algebraic subgroup of $\mathrm{SL}(2,\C)$.  Here, $G_1
=\mathcal{G}(k, \lambda)$ is the Galois group over $\C(z)$ of the
equation $L_2 = x'' - r (z) x = 0$. According to Theorem
\ref{abcritere}, we have to check if certain  specific primitive
integrals built from special function as Jacobi polynomials, are algebraic.

From Theorem~\ref{abcritere}, if $G_1$ is finite, then $G_2^{\circ}$ is
Abelian. In those cases, the existence of Jordan blocks with size $d=2$ does not give any obstacles for the integrability. This is why  we  are forced to look for such  obstructions considering Jordan blocks of size $d
= 3$. This problem is investigated  in
Section~4, where the results of Section~2 are also used. In this part of the paper we follow the
general ideas contained in Sections~2 and 3, but our considerations are much more technical.

In Section~5, we deal with the exceptional cases of potentials of degree $k =
\pm 2$, for which we prove that $G (\mathrm{VE}_t)^{\circ} \simeq G
(\mathrm{VE}_z)^{\circ}$ is Abelian. The strategy employed is completely
different and independent of the general frame of the paper. First, we  give a
direct proof of that result for $k = 2$. Then, we extract and discuss a general
principle of symmetry contained in Table \ref{tableMR}. Applying this principle,  we  deduce the following implication
\[
 G (\mathrm{VE}_t)^{\circ}\ \text{ Abelian for } k = 2 \quad \Longrightarrow\quad G
   (\mathrm{VE}_t)^{\circ}\ \text{is  Abelian for } k = - 2.
\]
For the non-expert reader, we should recommend to read this Section first,
since for $k = 2$, he shall see the frame of a very simple and particular VE.

In order to justify our study, in Section 6, we prove that the  Hessian matrix $V'' (c)$  for a homogeneous polynomial potential $V$ of degree $k$, can be an arbitrary symmetric matrix $A$ satisfying $Ac=(k-1)c$.
This is made by a dimensional arguments and study of complex symmetric
matrices.

\section{Theory for Jordan blocks of size two}

Let $(K, \partial)$ be an ordinary differential field with  constant subfield $\C$. We consider the following system of  two linear differential equations
over $K$.

\begin{align}
 x'' & =  rx,  \label{eqab1}\\
  y'' & = ry + sx.  \label{eqab2}
\end{align}

We  denote by $F_1$  and $F_2$   the   Picard-Vessiot
fields of equation~ \eqref{eqab1}, and the system \eqref{eqab1}-\eqref{eqab2}, respectively. The differential Galois group of extension $F_i/ K$ is denoted by $G_i$, for $i=1,2$.

We look for the conditions under  which $G_2^{\circ}$  is Abelian. Since $F_1$ may be seen as a subfield of
$F_2$, and $G_1$ as a quotient of $G_2$, we express these  conditions in
terms of $G_1^{\circ}$, and $r, s \in K$.

From now on, $\{x_1, x_2 \}$  denotes a basis of solutions of
\eqref{eqab1} normalised in such a way that
\begin{equation*}
 W (x_1 , x_2) = \det (X) = 1,  \mtext{ where}
 X =\begin{bmatrix}
      x_1 & x_2\\
     x_1' & x_2'
\end{bmatrix}.
\end{equation*}
For each $\sigma \in G_2$,  there exists matrix  $A (\sigma) \in \mathrm{SL}(2,\C)$, such that    $\sigma (X) = XA(\sigma)$.  Moreover, we chose $\{x_1, x_2 \}$  such
  that
\begin{itemize}
  \item if $G_1^{\circ} \simeq G_{\mathrm{a}}$,  then  for all
  $ \sigma \in G_1^{\circ}$, the matrix $A (\sigma)$ is a unipotent upper
  triangular matrix;

  \item if $G_1^{\circ} \simeq G_{\mathrm{m}}$,  then for all
   $ \sigma \in G_1^{\circ}$, the matrix $A (\sigma)$ is a diagonal matrix.
\end{itemize}
We recall here that in the above statements  $G_{\mathrm{a}}$ and $G_{\mathrm{m}}$ denote the additive and the multiplicative subgroups of $ \mathrm{SL}(2,\C)$.

 The group $G_1^{\circ}$ is a connected subgroup of  $\mathrm{SL}(2,\C)$. It is Abelian if  and only if it is  isomorphic either to $G_{\mathrm{a}}$, $G_{\mathrm{m}}$, or to $\{\id\}$.  Moreover, in \cite{kovacic86}  Kovacic proved the following
\begin{lemma}
  \label{abkovacic}
  \begin{enumerate}
    \item $G_1^{\circ} \simeq G_{\mathrm{a}}$, iff there exists a positive integer $m$ such
    that $x_1^m \in K$, and $x_2$ is transcendental over $K$. In this case the
    algebraic closure of $K$ in $F_1$ is $L = K [x_1]$.

    \item $G_1^{\circ} \simeq G_{\mathrm{m}}$, iff $x_1$ and $x_2$ are
    transcendental over $K$ but, $(x_1 x_2)^2 \in K$. In this case the
    algebraic closure of $K$ in $F_1$ is $L = K [x_1 x_2]$. Moreover, $L$ is at most
    quadratic over $K$.

    \item $G_1$ is a finite group if and only if both $x_1$ and $x_2$ are
    algebraic over $K$. Moreover, if this happens  then,
    $G_1$ is a finite subgroup of $\mathrm{SL} (2,\C)$ which is of one
    of the four types listed below:
\begin{enumerate}
\item Dihedral type: $G_1$ is conjugated to a finite subgroup of
\begin{equation*}
 D^{\dag} = \defset{
\begin{bmatrix}
    \lambda & 0\\
       0 & 1 / \lambda
\end{bmatrix}
}{\lambda\in\C^\star} \cup
 \defset{
\begin{bmatrix}
   0 & \lambda\\
       - 1 / \lambda & 0
\end{bmatrix}
}{ \lambda\in\C^\star} .
\end{equation*}
 \item
  Tetrahedral type: $G_1 /\{\pm \id\} \simeq \mathfrak{A}_4$.
\item
  Octahedral type: $G_1 /\{\pm \id\} \simeq \mathfrak{S}_4$.
\item
  Icosahedral type: $G_1 /\{\pm \id\} \simeq \mathfrak{A}_5$.

\end{enumerate}
\end{enumerate}
\end{lemma}
In the above,  $\mathfrak{S}_p$ and $\mathfrak{A}_p$ denote the symmetric, and the alternating group of $p$ elements, respectively.
\begin{definition}
  \label{abcondition}Let  $\varphi = \int sx_1^2$ and $\psi =
  \int x_1^{-2}$. We  define  the following conditions
  \begin{description}
    \item[$(\alpha)$]There exists $c \in \C$ such that $\varphi + c \psi \in
    L$.

    \item[$(\beta)$] $\varphi  \psi - 2 \int \varphi \cdot \psi' = 2 \int
    \varphi'  \psi - \varphi \cdot \psi \in L [\psi]$.

    \item[$(\gamma)$] There exists $\phi_1 \in L$ such that $(\phi_1 x_1^2)' =
    sx_1^2$.

    \item[$(\delta)$] There exists $\phi_2 \in L$ such that $(\phi_2 x_2^2)' =
    sx_2^2$.
  \end{description}
\end{definition}
With the above  notations and definitions  our main result in this section is the following.
\begin{theorem}
  \label{abcritere} The group $G_2^{\circ}$ is  Abelian if and only if  one of the following cases occur
  \begin{enumerate}
    \item $G_1$ is a finite group.

    \item $G_1^{\circ} \simeq G_{\mathrm{a}}$ and condition  $(\beta)$
    holds.

    \item $G_1^{\circ} \simeq G_{\mathrm{m}}$ and conditions $(\gamma)$ and $(\delta)$
    hold.
  \end{enumerate}
 When $G_1^{\circ} \simeq G_{\mathrm{a}}$,  then   $(\beta) \Rightarrow (\alpha)$.  Hence,  $ (\alpha)$ is a necessary condition for  $G_2^{\circ}$ to be Abelian in this case.
\end{theorem}

The proof of the theorem will be done at the end of this section. Before, for
sake of clarity, we explain the main ideas of the proof.

A good illustration of  the Kronecker observation in arithmetic is the following
example. Let
\[ f (X) = X^3 + pX + q \in \Q[X], \]
be an irreducible polynomial. Let $\mathrm{Gal} (f /\Q)$ be its Galois
group over the rationals, and $\Delta = - 4 p^3 - 27 q^2$ be the discriminant
of $f$. The group $\mathrm{Gal} (f /\Q)$ can be either $\mathfrak{S}_3$, or
$\mathfrak{A}_3$.
Moreover,
\[ \mathrm{Gal} (f /\Q) \simeq \mathfrak{A}_3 \ \Longleftrightarrow\  \Delta
   \in (\Q)^2 . \]
In other words, $\mathrm{Gal} (f /\Q)$ is Abelian, iff $\Delta$ is a
square of a rational number.  This is,  in the considered example,  the precise ``arithmetical
condition'' that governs the Abelianity of the Galois group.

In differential Galois theory, the analogue of the discriminant is the
Wronskian determinant. Therefore our idea was to express the Abelianity
condition for $G_2^{\circ}$ in terms of certain properties of the Wronskian
determinant.

We proceed in the three following steps
\begin{enumerate}
  \item The very specific form of system \eqref{eqab1}--\eqref{eqab2}, allows
   to express the Abelianity of $G_2^{\circ}$ in terms of its subgroup $H =
  \mathrm{Gal}_{\partial} (F_2 / F_1)$.

  \item Next,  we translate  this  group conditions  into properties of certain Wronskians.

  \item In a third step, thanks to Lemma \ref{ostrowski}, we express these
  Wronskian properties in terms of the algebraicity of certain primitive integrals.
\end{enumerate}
Later, in the applications, we shall not use the case of Theorem
\ref{abcritere} where $G_1^{\circ} \simeq G_{\mathrm{m}}$. This is because those cases
only happen for potentials of degree $k = \pm 2$ for which other kind of
arguments will be applied in Section 5. Therefore, at first reading, this part
of the proof of Theorem \ref{abcritere} may be avoided.

\subsection{Group formulation of the criterion}

The system \eqref{eqab1} and \eqref{eqab2} may be written into the  matrix form:
\begin{equation}
\label{j2mat}
 \begin{bmatrix}
      x'\\
     x''\\
     y'\\
     y''
 \end{bmatrix}=
\begin{bmatrix}
    0 & 1 & 0 & 0\\
     r & 0 & 0 & 0\\
     0 & 0 & 0 & 1\\
     s & 0 & r & 0
\end{bmatrix}
\begin{bmatrix}
     x\\
     x'\\
     y\\
     y'
\end{bmatrix}
=
\begin{bmatrix}
   R & 0\\
     S & R
\end{bmatrix}
\begin{bmatrix}
     x\\
     x'\\
     y\\
     y'
\end{bmatrix},
\end{equation}
where
\begin{equation}
 \label{eq:RS}
R:=\begin{bmatrix}
 0 & 1\\
r & 0
\end{bmatrix}, \quad
S:=\begin{bmatrix}
 0 & 0\\
s& 0
\end{bmatrix}.
\end{equation}
For a given basis $\{x_1 , x_2 \}$ of  solutions of equation~\eqref{eqab1}, we set
\begin{equation*}
 X=\begin{bmatrix}
 x_1 & x_2\\
  x_1' & x_2'
\end{bmatrix} \mtext{and}
 Y=\begin{bmatrix}
 y_1 & y_2\\
  y_1' & y_2'
\end{bmatrix},
\end{equation*}
where $y_1$ and $y_2$  are two particular solutions of
\eqref{eqab2}, that is:
\begin{equation}
\label{y12}
\left\{
\begin{split}
      y_1'' &= ry_1 + sx_1\\
     y_2'' &= ry_2 + sx_2
\end{split}
\right.
\end{equation}
Then the following $4 \times 4$ matrix
\begin{equation*}
 \Xi_2=
\begin{bmatrix}
     X & 0\\
     Y & X
\end{bmatrix},
\end{equation*}
 is a fundamental matrix of solutions of~\eqref{j2mat}.
For each $ \sigma \in G_2$,   we have
\begin{equation}
\label{ms}
 \sigma (\Xi_2) =
\begin{bmatrix}
   \sigma (X) & 0\\
  \sigma (Y) & \sigma (X)
\end{bmatrix}=
\Xi_2 M(\sigma).
\end{equation}
Performing the above multiplication we can easily notice that
the $4 \times 4$ matrix $M (\sigma)$ has the  form
\begin{equation*}
  M (\sigma)=
\begin{bmatrix}
   A (\sigma) & 0\\
  B (\sigma) & A (\sigma)
\end{bmatrix}.
\end{equation*}
Therefore, $G_2$ can be identified with a subgroup of $\mathrm{SL}(4,\C)$:
\begin{equation}
 \label{eqab4}
  G_2 \subset G_{\max} =\defset{%
 M(A,B)= \begin{bmatrix}
      A & 0\\
    B & A
\end{bmatrix}
}{%
 A \in \mathrm{SL}(2,\C), B \in \M(2,\C) }.
\end{equation}
For   $M_i=M_i(A_i, B_i)\in G_{\max}$, with $i=1,2$, we have
\begin{equation}
 \label{eqab5}
 M_1 M_2 = \begin{bmatrix}
    A_1 A_2 & 0\\
    B_1 A_2 + A_1 B_2 & A_1 A_2
 \end{bmatrix}.
\end{equation}
\begin{definition}
  \label{abH}We denote  by $H := \mathrm{Gal}_{\partial} (F_2 / F_1)$.  The groups  $H_{\max}$, $H_{\mathrm{a}}$ and $H_{\mathrm{m}}$
are subgroups of $G_{\max}$, defined by,
  \begin{equation*}
   H_{\max}:=\defset{%
N (B) =\begin{bmatrix}
           \id & 0\\
      B & \id
       \end{bmatrix}
}{%
 B \in \M(2,\C)},
  \end{equation*}
\begin{equation*}
 H_{\mathrm{a}} :=\defset{ N(B)\in H_{\max} }{%
B = \begin{bmatrix}
     a & b\\
        0 & a
\end{bmatrix}, \ a, b\in \C },
\end{equation*}
\begin{equation*}
 H_{\mathrm{m}} :=\defset{ N(B)\in H_{\max} }{%
B = \begin{bmatrix}
     a & 0\\
        0 & d
\end{bmatrix}, \ a, d\in \C }.
\end{equation*}
\end{definition}
From \eqref{eqab5}, we have $N (B_1) N (B_2) = N (B_1 + B_2)$. So, $H_{\max}$ is a
vector group of dimension 4 isomorphic to $(\M(2,\C), +)$.
\begin{proposition}
  \label{abG}With the  notations above we have the following.
  \begin{enumerate}
    \item The Picard Vessiot extension $F_2 / F_1$ is a regular fields extension  and its Galois
    group $H = \mathrm{Gal}_{\partial} (F_2 / F_1)$ is a vector group.

    \item The algebraic closure of $K$ in $F_2$ coincides with the
    algebraic closure  $L$ of $K$ in $F_1$.

    \item The kernel of the restriction map $\mathrm{Res}^{\circ} : G_2^{\circ}
    \rightarrow G^{\circ}_1$ coincides with $H = \mathrm{Gal}_{\partial} (F_2 /
    F_1)$.

    \item If $G_1$ is finite then $G_2^{\circ} = H = \mathrm{Gal}_{\partial}
    (F_2 / F_1)$ is Abelian.
  \end{enumerate}
\end{proposition}
\begin{proof}
(1)~Let $\mathrm{Res} : G_2 \rightarrow G_1$,  $\sigma \mapsto
  \sigma |_{F_1}$, be the restriction map.  We have
\begin{equation*}
 M (\sigma) =
\begin{bmatrix}
    A (\sigma) & 0\\
    B (\sigma) & A (\sigma)
\end{bmatrix}\longmapsto  \mathrm{Res}(M(\sigma))=  A (\sigma)
\end{equation*}
Since  $H
  = \mathrm{Gal}_{\partial} (F_2 / F_1) = \mathrm{Ker} (\mathrm{Res})$, the algebraic subgroup $H$ of $G_2$  may be viewed
  as an algebraic subgroup of $H_{\max}$. It is therefore a vector group,
  hence connected.

(2)~Let $u \in F_2$ be algebraic over $K$. Since $H$ is connected, $Hu
  =\{u\}$, and thus $u \in F_1$ is algebraic over $K$.

(3)~Since the restriction map $\mathrm{Res} : G_2 \rightarrow G_1$ is a
  surjective morphism of algebraic groups, it maps $G_2^{\circ}$ onto
  $G^{\circ}_1$. Denoting  by $\mathrm{Res}^{\circ}$ the restriction of $ \mathrm{Res}$ to $G_2^\circ$ and
putting $H' = \mathrm{Ker} (\mathrm{Res}^{\circ})$,   we have the following
  commutative diagram of algebraic groups, whose lines are exact sequences
\begin{equation*}
 \begin{CD}
     \{\id\} @>>> H @>>> G_2 @>\mathrm{Res}>> G_1 @>>> \{\id\} \\
   @.       @AAA       @AAA                           @AAA      @.  \\
  \{\id\} @>>> H' @>>> G_2 ^\circ @>\mathrm{Res^\circ}>> G_1 ^\circ@>>> \{\id\} \\
 @.       @AAA       @AAA                           @AAA      @.  \\
  \{\id\} @>>> \{\id\} @>>> \{\id\} @>>> \{\id\} @>>> \{\id\}
\end{CD}
\end{equation*}
  Applying the snake lemma to the  first two lines we obtain following exact
  sequence
\begin{equation*}
\begin{CD}
   \{\id\}@>>>H / H' @>>> G_2 / G_2^{\circ} @>>>G_1 / G_1^{\circ} @>>> \{\id\}
     \end{CD}
\end{equation*}
 But $G_2 / G_2^{\circ}$ is finite, so  $H / H'$ is also finite. Moreover, $H / H'$ as a quotient
  of vector group is also a vector group hence, it is the trivial vector
  group. That is $H' = H$, and $G_2 / G_2^{\circ}$ is isomorphic to $G_1 /
  G_1^{\circ}$. Moreover, the second line of the commutative diagram reduces
  to the exact sequence
  \begin{equation}
    \begin{CD}
      \{\id\} @>>> H = \mathrm{Gal}_{\partial} (F_2 / F_1) @>>>
      G_2^{\circ}@>>>  G^{\circ}_1 @>>> \{\id\}.
    \end{CD} \label{eqab6}
  \end{equation}
 (4)~If $G_1$ is finite, $G_1^{\circ} =\{\id \}$ in \eqref{eqab6}, so
  $G_2^{\circ} = H$ is Abelian.
\end{proof}
In general, from~\eqref{eqab6} we have
\begin{equation*}
 G_2^{\circ} \subset \defset{M(A,B)\in G_{\max}}{A \in G_1^{\circ}},
\end{equation*}
 and, if $G_2^\circ$  is
Abelian, then $G_1^{\circ}$ is an Abelian algebraic subgroup of
$\mathrm{SL}(2,\C)$. If this happens, then $G_1^{\circ}$ is  isomorphic either to
$\{\id\}$, or $G_{\mathrm{a}}$,  or $G_{\mathrm{m}}$.

If $G_1^{\circ} =\{\id\}$ we have
seen above that $G_2^{\circ}$ is an Abelian vector group.
However, if  $G_1^{\circ}$ is  isomorphic  either to $G_{\mathrm{a}}$ or $G_{\mathrm{m}}$, then we have to
find  conditions under which $G_2^{\circ}$ is  Abelian. For
that purpose we need the following conjugation formula, which is obtained from
\eqref{eqab5}  by direct computations. Namely, for all  $ M(A,B)\in G_{\max}$  and  all  $N(C) \in H_{\max}$,
we have
\begin{equation}
  M(A,B)N(C)M(A,B)^{- 1} = N (ACA^{-1}).  \label{eqab7}
\end{equation}

\begin{proposition}
  \label{abcritereH}
If $G_1^{\circ}$ is  isomorphic either to
  $G_{\mathrm{a}}$ or $G_{\mathrm{m}}$, then $G_2^{\circ}$ is Abelian if and only if $H =
  \mathrm{Gal}_{\partial} (F_2 / F_1)\subset Z (G_2^{\circ}) $, i.e., $H$ is contained in the center of
  $G_2^{\circ}$.
  \begin{enumerate}
    \item If $G_1^{\circ} \simeq G_{\mathrm{a}}$, then $G_2^{\circ}$ is Abelian iff $H$ is
    a subgroup of $H_{\mathrm{a}}$.

    \item If $G_1^{\circ} \simeq G_{\mathrm{m}}$, then $G_2^{\circ}$ is Abelian iff $H$ is
    a subgroup of $H_{\mathrm{m}}$.
  \end{enumerate}
Here  $H_{\mathrm{a}}$  and  $H_{\mathrm{m}}$ are the groups defined in Definition~\ref{abH}.
\end{proposition}

\begin{proof}
  If $G_2^{\circ}$ is Abelian, then $H \subset Z (G_2^{\circ})$. Conversely,
  let us assume that $G_1^{\circ}$ is   isomorphic either to $G_{\mathrm{a}}$, or $G_{\mathrm{m}}$, and $H
  \subset Z (G_2^{\circ})$. For any $M_0 \in G_2^{\circ} \backslash H$, the
  subgroup $\Omega$ generated by $M_0$ and $H$,  as well as its Zariski closure $\bar{\Omega}$, is an Abelian subgroup of $G_2^{\circ}$. By equation \eqref{eqab6}, we have $\dim G_2^{\circ} = \dim H + 1$, but
  \[ \left( \dim H + 1 \leq \dim \bar{\Omega}
     \leq \dim G_2^{\circ} \right) \quad \Longrightarrow\quad \dim \bar{\Omega} = \dim
     G_2^{\circ} . \]
  Since $G_2^{\circ}$ is connected, we deduce that  $G^{\circ}_2 = \bar{\Omega}$ is Abelian.

Note that $H \subset Z (G_2^{\circ})$ iff for all $M=M(A,B)  \in G_2^{\circ}$, and all $N=N(C)  \in H$, we have, thanks to~\eqref{eqab7},
\begin{equation*}
 MNM^{- 1} = N   \Longleftrightarrow N(ACA^{-1})=N(C) \Longleftrightarrow ACA^{-1}=C \Longleftrightarrow [A,C]=0.
\end{equation*}
Now, we can prove the remaining points.

(1)~If $G_1^{\circ} \simeq G_{\mathrm{a}}$ we put
\begin{equation*}
 A = A (t) =
\begin{bmatrix}
     1 & t\\
    0 & 1
\end{bmatrix}
\mtext{and} C=
\begin{bmatrix}
        a & b\\
    c & d
\end{bmatrix}.
\end{equation*}
We have
\begin{equation*}
  [C, A] =
\begin{bmatrix}
            - tc & t (a - d)\\
       0 & - tc
\end{bmatrix},
\end{equation*}
and thus
\begin{equation*}
 ( \forall\; A \in G_{\mathrm{a}}\ \ [C , A] = 0)\  \Longleftrightarrow\  (c = 0 \mtext{and} a = d)\
     \Longleftrightarrow C =
\begin{bmatrix}
        a & b\\
    0 & a
\end{bmatrix}\ \Longleftrightarrow\  H\subset H_{\mathrm{a}}.
\end{equation*}

(2)~If $G_1^{\circ} \simeq G_{\mathrm{m}}$, we proceed in a similar way and  we  put
\begin{equation*}
 A = A (t) =
\begin{bmatrix}
     t & 0\\
    0 & 1/t
\end{bmatrix}
\mtext{and} C=
\begin{bmatrix}
        a & b\\
    c & d
\end{bmatrix}.
\end{equation*}
Then we obtain
\begin{equation*}
  [C, A] =
\begin{bmatrix}
             0 & - b\left (t - {1}/{t}\right)\\
       c \left(t - {1}/{t}\right) & 0
\end{bmatrix},
\end{equation*}
and thus
\begin{equation*}
 ( \forall\; A \in G_{\mathrm{m}} \ \ [C,  A] = 0) \Longleftrightarrow (b=0 \mtext{and} c = 0)
     \Longleftrightarrow C =
\begin{bmatrix}
        a & 0\\
    0 & d
\end{bmatrix} \Longleftrightarrow H \subset H_{\mathrm{m}}.
\end{equation*}
\end{proof}

\subsection{From group to Wronskian relations}

For two elements $f, g \in F_2$, we set
\begin{equation*}
 W (f , g) = \begin{vmatrix}
  f & g\\
  f' & g'
\end{vmatrix}.
\end{equation*}
Observe that for $\sigma \in G_2$ we have
\begin{equation*}
\sigma ( W (f , g)) = \begin{vmatrix}
  \sigma(f) & \sigma(g)\\
  \sigma(f)' & \sigma(g)'
\end{vmatrix}=W(  \sigma(f) , \sigma(g)).
\end{equation*}
Let  $x_1$, $x_2$, $y_1$, $y_2$, $X$  and $Y$  be as  they were defined in Section 2.1.
\begin{definition}
  \label{defW}We define the following three conditions
  \begin{description}
    \item[$W_1:$] \  $W (x_1 , y_1) \in F_1$.

    \item [$W_2:$]\ $W (x_2 , y_2) \in F_1$.

    \item[$W_3:$] \  $ x_1 W (x_1 , y_2) - y_1 \in F_1$.
  \end{description}
\end{definition}

\begin{proposition}
  \label{abcritereW}Let $H_{\mathrm{a}}$ and $H_{\mathrm{m}}$ be the groups given in
  Definition \ref{abH}.  We have
  \begin{enumerate}
    \item $H = \mathrm{Gal}_{\partial} (F_2 / F_1) \subset H_{\mathrm{a}}$ iff
    condition $W_3$ is fulfilled. Moreover, $W_3 \Rightarrow W_1$.

    \item $H = \mathrm{Gal}_{\partial} (F_2 / F_1) \subset H_{\mathrm{m}}$ iff the
    conditions $W_1$ and $W_2$ are fulfilled.
  \end{enumerate}
\end{proposition}

\begin{proof}
  Let $\sigma\in H$.  Then $\sigma(X)=X$,  and
  $\sigma (Y) = Y + XB (\sigma)$ for a certain
\begin{equation*}
 B (\sigma) =
\begin{bmatrix}
   a & b\\
    c & d
\end{bmatrix}
 \in \M(2,\C).
\end{equation*}
 Therefore, the action of $\sigma$ on $Y$
  is given by the relations
 \begin{equation*}
\left\{
  \begin{split}
   \sigma (y_1) & = & y_1 + ax_1 + cx_2,\\
       \sigma (y_2) & = & y_2 + bx_1 + dx_2.
\end{split}\quad \right.
 \end{equation*}
  From these relations, the action of $\sigma$ on the Wronskians is given by the following
  formulae
\begin{align*}
      \sigma \left( W (x_1 , y_1) \right) & = \begin{vmatrix}
         x_1 & y_1 + ax_1 + cx_2\\
         x_1' & y_1' + ax_1' + cx_2'
       \end{vmatrix} = W (x_1 , y_1) + c,\\
       \sigma \left( W (x_2 , y_2) \right) & =  \begin{vmatrix}
         x_2 & y_2 + bx_1 + dx_2\\
         x_2' & y_2' + bx_1' + dx_2'
       \end{vmatrix}  = W (x_2 , y_2) - b,\\
       \sigma \left( W (x_1 , y_2) \right) & = \begin{vmatrix}
         x_1 & y_2 + bx_1 + dx_2\\
         x_1' & y_2' + bx_1' + dx_2'
       \end{vmatrix}  =  W (x_1 , y_2) + d,\\
       \sigma \left( W (x_2 , y_1) \right) & = \begin{vmatrix}
         x_2 & y_1 + ax_1 + cx_2\\
         x_2' & y_1' + ax_1' + cx_2'
       \end{vmatrix}  =  W (x_2 , y_1) - a.
\end{align*}
  To obtain these formulae, we  used
  the fact that $W (x_1 , x_2) = 1$.
 Moreover, we   also  have
  \begin{eqnarray*}
    \sigma \left( x_1 W (x_1 , y_2) - y_1) \right. & = & x_1 \sigma (W (x_1 ,
    y_2)) - \sigma (y_1)\\
    & = & x_1 W (x_1 , y_2) + dx_1 - (y_1 + ax_1 + cx_2),\\
    \sigma \left( x_1 W (x_1 , y_2) - y_1) \right. & = & [x_1 W (x_1 , y_2) -
    y_1] + (d - a) x_1 - cx_2 .
  \end{eqnarray*}
  Therefore,
\begin{gather*}
  \sigma \left( W (x_1 , y_1) \right)  = W (x_1 , y_1)
 \     \Longleftrightarrow \  c = 0, \\
 \sigma \left( W (x_2 , y_2) \right) =  W (x_2 , y_2) \
     \Longleftrightarrow\  b = 0, \\
   \sigma \left( x_1 W (x_1 , y_2) - y_1) \right. =   x_1 W (x_1 , y_2)
       - y_1\ \Longleftrightarrow\  [d = a \mtext{and}c = 0].
\end{gather*}
  For the last equivalence we used  the fact that $x_1$ and $ x_2 $ are
  $\C$-linearly independent.

  From the above equivalences we deduce that for  $ \sigma \in H$,
we have
\begin{gather*}
 \sigma \in H_{\mathrm{m}} \   \Longleftrightarrow \   \left( \sigma ( W (x_1 , y_1) )
       = W (x_1 , y_1)  \mtext{and}\sigma ( W (x_2 , y_2) ) = W (x_2 ,
       y_2)\right),\\
\sigma \in H_{\mathrm{a}} \ \Longleftrightarrow \sigma ( x_1 W (x_1 , y_2) -
       y_1) = x_1 W (x_1 , y_2) - y_1.
\end{gather*}
and,
  moreover,
  \[ \left( \sigma ( x_1 W (x_1 , y_2) - y_1) = x_1 W (x_1 , y_2) - y_1\right) \
\Longrightarrow \  \left( \sigma ( W (x_1 , y_1) ) = W (x_1 , y_1)\right) .
  \]
  But  for $f \in F_2$ we have
\[
 f \in F_1\quad \Longleftrightarrow \quad  \left( \forall \sigma \in H,
  \sigma (f) = f\right).
\]
So, $H \subset H_{\mathrm{m}}$ if and only if conditions $W_1$ and
  $W_2$ hold. Similarly,  $H \subset H_{\mathrm{a}}$ if and only if the conditions $W_3$
  is  satisfied, moreover   $W_3 \Rightarrow W_1$.
\end{proof}

\subsection{From Wronskian to integral relations}

\subsubsection{Computation of the Wronskian and resolution of equation~\eqref{eqab2}}

From Proposition \ref{abG} we know that $H = \mathrm{Gal}_{\partial} (F_2 / F_1)$ is
a vector group, so it is solvable. Therefore, by Liouville-Kolchin solvability
theorem, equation \eqref{eqab2}: $y'' = ry + sx$ can be solved by finite
integrations.  In the following lemma   we give, among others things, the explicit form of its solutions.
Notice that in this lemma we do not make any assumption about the group $G_1$.

Let $x_1$ be a non-zero solution of equation~\eqref{eqab1}. According to Definition~\ref{abcondition}, let us set
$\varphi = \int sx_1^2$ and $\psi =
  \int x_1^{-2}$.  Then $x_2=x_1\psi$ is another solution of~\eqref{eqab1}, and $W(x_1,x_2)=1$.  Let $y_1$  and $y_2$ be two particular solution of~\eqref{eqab2} given~\eqref{y12}. Then we have the following.
\begin{lemma}
  \label{abcomputation}Up to additive constants we have
\begin{equation*}
 \begin{split}
 W (x_1 , y_1) &= \int sx_1^2,   \qquad\ \ \  W (x_2 , y_2) = \int sx_2^2 , \\
W (x_1 , y_2)&= \int sx_1 x_2 , \qquad  W (x_2 , y_1) = \int sx_1 x_2,
\end{split}
\end{equation*}
\begin{equation*}
y_1 =  x_1  \int \left( x_1^{-2} \int sx_1^2\right) ,\qquad
   y_2  = x_2  \int
     \left(x_2^{-2} \int sx_2^2\right),
\end{equation*}
and
  \[
  x_1 W (x_1 , y_2) - y_1 = x_1 Q, \mtext{where}  Q =  \varphi  \psi - 2 \int \varphi  \psi' =  - \varphi  \psi +2  \int \psi\varphi'.
\]
\end{lemma}
\begin{proof}
Identities with Wronskians can be checked by a direct differentiation. Formulae for $y_1$ and $y_2$ are obtained by a classical variations of constants method.
\end{proof}
\begin{corollary}
  \label{abtrace}
Let $\sigma\in H = \mathrm{Gal}_{\partial} (F_2 / F_1)$, and $N(B(\sigma))\in H_{\max}$ be the matrix of $\sigma$.  Then $\mathrm{Tr} (B(\sigma)) = 0$.
\end{corollary}
\begin{proof}
As in the proof of  Proposition \ref{abcritereW} we set $B=B(\sigma)=\bigl(\begin{smallmatrix} a &b \\ c &d \end{smallmatrix}\bigr)$. From this proof we know also that
  \[ \sigma \left( W (x_1 , y_2) \right) = W (x_1 , y_2) + d \mtext{ and }
     \sigma \left( W (x_2 , y_1) \right) = W (x_2 , y_1) - a. \]
  But from Lemma \ref{abcomputation} we know that $\Delta = W (x_1 , y_2) - W (x_2 , y_1)$
  is a constant belonging to $\C$. Therefore, for  $\sigma \in H$,
  \[ \sigma (\Delta) = \Delta + d + a = \Delta + \mathrm{Tr} (B) = \Delta. \]
  So, $\mathrm{Tr} (B) = 0$.
\end{proof}

\subsubsection{Study of the conditions $W_i$}

\begin{lemma}
  \label{abcritereint}
   Let $T (F_1 / K) \subset F_1$ be the Picard Vessiot ring of $F_1 /
    K$. With the notations of Lemmas \ref{ostrowski} and \ref{abkovacic}, we
    have
\begin{enumerate}
\item     \begin{itemize}
      \item If \ $G_1^{\circ} \simeq G_{\mathrm{a}}$, then $T (F_1 / K) = L [x_2] = L
      [\psi]$.

      \item If  $G_1^{\circ} \simeq G_{\mathrm{m}}$, then $T (F_1 / K) = L [x_1, x_2] = L
      [x_1, x_1^{-1}]$.
    \end{itemize}
    \item If $G_1^{\circ} \simeq G_{\mathrm{a}}$, then the condition $W_1$ is equivalent to
    $(\alpha)$ and $W_3$ to $(\beta)$.

    \item If \ $G_1^{\circ} \simeq G_{\mathrm{m}}$, then the condition $W_1$ is equivalent to
    $(\gamma)$ and $W_2$ to $(\delta)$.
  \end{enumerate}
\end{lemma}

\begin{proof}
  (1)~From now on, as in Lemma \ref{abkovacic}, $L$  denotes the algebraic
  closure of $K$ in $F_1$. From the relation $W (x_1 , x_2) = 1$, and Lemma
  \ref{ostrowski}, we have
  \[ T (F_1 / K) = K [x_1, x_2, x_1', x_2'] = L [x_1, x_2, x_1', x_2'] . \]
  Moreover, $F_1$ is the field of fractions of the ring $T (F_1 / K)$. Let us
  compute this ring in the two particular cases.

  If $G_1^{\circ} \simeq G_{\mathrm{a}}$, then, by assumption, $x_1 \in L$, so $x_1' \in L$.
  Since $x_2 =  x_1 \int x_1^{-2} =x_1\psi$, we have
  \[ T (F_1 / K) = L [x_2] = L [  x_1\psi] =L[\psi]. \]

If $G_1^{\circ} \simeq G_{\mathrm{m}}$, then $G_1^{\circ}$ acts on
  $x_1$ by a character, so  it acts on $x_1'$ by the same character. Therefore, the
  logarithmic derivative $x_1' / x_1$ is left invariant by $G_1^{\circ}$,
  hence belongs to $L$. Moreover, from Lemma \ref{abkovacic}, $x_1 x_2 \in L$,
  so similarly $x_2' / x_2 \in L$, and we have
  \[ T (F_1 / K) = L [x_1, x_2] = L [x_1, {x_1}^{-1}] . \]

 (2)~Since $G_1^{\circ} \simeq G_{\mathrm{a}}$,
  \[ F_1 = L (x_2) = L ( \int x_1^{-2}) = L (\psi) . \]
  Therefore, from Lemma \ref{abcomputation}, the condition $W_1$ may be written
  \[ W (x_1 , y_1) = \int sx_1^2 \in L ( \int x_1^{-2}) . \]
  Thus, condition $W_1$ implies that the two primitive integrals: $\varphi = \int sx_1^2$,
  and $\psi = \int x_1^{-2}$ are algebraically dependant over $L$. Hence, by the
   Ostrowski-Kolchin theorem (Lemma \ref{ostrowski} point 1), this implies that there
  exists $(c_1, c_2) \in \C^2 \backslash \{0, 0\}$ such that
  \[ c_1 \varphi + c_2 \psi \in L. \]
  But $c_1 = 0$ implies that $\psi = \int x_1^{-2} \in L$, and $x_2 = - x_1
  \psi$ is algebraic over $K$, however it is not true. So, dividing the linear
  relation by $c_1$ we get that $W_1 \Rightarrow (\alpha)$. Conversely, if
  $(\alpha)$ holds then, $\int sx_1^2 \in L ( \int x_1^{-2}) = L (x_2)
  = F_1$ and $W_1$ is satisfied.

  From Lemma \ref{abcomputation} we have
  \begin{equation*}
    x_1 Q  =  x_1 W (x_1 , y_2) - y_1,
  \end{equation*}
where
\begin{equation*}
 Q  =  \varphi  \psi - 2 \int \varphi  \psi' =  - \varphi  \psi +2  \int \psi\varphi'.
\end{equation*}

  From Lemma \ref{ostrowski}, the element $x_1 Q = x_1 W (x_1 , y_2) - y_1 \in
  T (F_2 / K)$. But here, $x_1$ is algebraic over $K$ so from Lemma
  \ref{ostrowski} again, $Q = \frac{1}{x_1} \cdot x_1 Q \in T (F_2 / K)$.
  Therefore, the condition $Q \in F_1$, is equivalent to $Q \in T (F_1 / K)$,
  because $T (F_1 / K)$ is the algebra containing the elements of $F_1$ which
  are solutions of a certain linear differential equation over $K$. So, we have the following
  equivalences
  \[
W_3 \ \Longleftrightarrow\ ( x_1 Q \in F_1)\  \Longleftrightarrow \ (Q \in F_1)
    \ \Longleftrightarrow\ (Q \in T (F_1 / K)) \ \Longleftrightarrow\  (Q \in L [\psi]) . \]
 Thus, condition  $W_3$ is equivalent to condition $(\beta)$.

 (3)~Since $G_1^{\circ} \simeq G_{\mathrm{m}}$, the role of $x_1$ and $x_2$ are
  symmetric. We have to prove only  that  conditions $W_1$ and $(\gamma)$ are equivalent. As before, $W (x_1 ,
  y_1) \in T (F_2 / K)$ so,
  \[
W (x_1 , y_1) \in F_1 \ \Longleftrightarrow\  W (x_1 , y_1) \in T (F_1 / K) = L
     [x_1 , x_1^{-1}] . \]
  Since $W (x_1 , y_1) = \int sx_1^2$, condition $W_1$ is equivalent to
  \[ \int sx_1^2 \in L [x_1 ,x_1^{-1}] . \]
  The above condition is fulfilled iff we have a relation of the form
  \[ \int sx_1^2 = \sum_{n=p}^q f_n x_1^n,  \quad p\leq q; \ p,q\in\Z, \]
  with $f_n \in L$.  Differentiating  the above equation we obtain
  \[ sx_1^2 = \sum_{n=p}^q (f_n' + f_n n \theta) x_1^n , \]
where $\theta=x_1' / x_1 \in L$.

  But $x_1$ is transcendental over $L$,  so from the last  formula we have
  $f_n' + f_n n \theta = 0$  for $ n \neq 2$, and $s = f_2' + 2 f_2 \theta$.
  Thus we have
  \[ sx_1^2 = (\phi_1 x_1^2)', \]
  with $\phi_1 = f_2 \in L$. This proves that condition $W_1$ is equivalent to
 condition  ($\gamma$).
\end{proof}

\subsection{Proof of Theorem \ref{abcritere}}

\begin{proof}
 As a connected subgroup of $\mathrm{SL}(2 ,\C)$, group $G_1^{\circ}$ is
  isomorphic to one of the following groups: $\{\id\}, G_{\mathrm{a}}, G_{\mathrm{m}}$, the semi-direct product $G_{\mathrm{m}}\ltimes G_{\mathrm{a}}
  $ or $\mathrm{SL}(2 ,\C)$. If $G_1^{\circ}$ is  Abelian, then the last two
  possibilities must be excluded.
  \begin{itemize}
    \item If $G_1^{\circ} =\{\id\}$, then $G_2^{\circ}$ is Abelian thanks
    to point 5 of Proposition \ref{abG}.

    \item If $G_1^{\circ} \simeq G_{\mathrm{a}}$, (resp. $G_1^{\circ} \simeq G_{\mathrm{m}}$), then
    the proof follows from Proposition \ref{abcritereW} and point 2, (resp.
    point 3) of Lemma \ref{abcritereint}.
  \end{itemize}
\end{proof}

\section{Elimination of the Jordan blocks with $G_1^{\circ} \simeq G_{\mathrm{a}}$.}

We now apply the results of the previous section to the study of the
connected component $G (\mathrm{VE}_z)^{\circ}$ of the Galois Galois group of
the VE (\ref{eqVEredz})
\[ \frac{d^2 \zeta}{dz^2} = [ r_0 (z) \id + s (z) V'' (c)] \zeta . \]
Our main result in this section is the following.
\begin{theorem}
  \label{jordanGa}Assume that $V'' (c)$ has a Jordan block of size
  $d \geq 2$, and $G_1^{\circ} \simeq G_{\mathrm{a}}$.  Then $G (\mathrm{VE}_z)^{\circ}$
  is not Abelian. This corresponds to the elimination of rows 2,3, and 4 in
  Table \ref{tableMR}.
\end{theorem}

\begin{remark}
  \label{jordanredGa}Let $B (\lambda, d)$ be a Jordan block of
  $V'' (c)$ with size $d \geq 2$ and eigenvalue $\lambda$. Since
  $G_1^{\circ}$ is isomorphic to $G_{\mathrm{a}}$ and corresponds to the VE
  \[ \frac{d^2 \eta}{dt^2} = - \lambda \varphi^{k - 2} (t) \eta, \]
  we deduce from Theorem \ref{MR} that necessarily, the pair $(k, \lambda)$
  must belong to row 2, 3, or 4 in Table \ref{tableMR}. Now, passing to the VE in
  the $z$ variable, we know that the system \eqref{eqblock} with Galois group
  $G_2$
  \[ \begin{bmatrix}
       x''\\
       y''
     \end{bmatrix} = \begin{bmatrix}
       r & 0\\
       s & r
     \end{bmatrix}\begin{bmatrix}
       x\\
       y
     \end{bmatrix}, \]
  is a subsystem of VE (\ref{eqVEredz}). From Proposition \ref{VEtVEz}, $G
  (\mathrm{VE}_t)^{\circ} \simeq G (\mathrm{VE}_z)^{\circ}$ and $G_2^{\circ}$ is a
  quotient of $G (\mathrm{VE}_z)^{\circ}$. Therefore, it is enough   to prove
  that $G_2^{\circ}$ is not Abelian. To this aim we proceed as follows. According to Theorem~\ref{abcritere}, we have to prove that condition $(\beta)$ is not fulfilled.  Since  $(\beta)\Rightarrow(\alpha)$,  and  $(\alpha)$ is much easier to check than $(\beta)$, at first we check if $(\alpha)$ is  fulfilled. Since $(\alpha)$ is a  condition concerning the  primitive integrals
$\varphi=\int s x_1^2$ and $\psi=\int x_1^{-2}$, where $x_1$ is an algebraic solution of equation $x''=rx$, we first have to investigate  analytical  properties of these integrals.
\end{remark}

\subsection{Assumptions and notations}

We assume that $G_1^{\circ} \simeq G_{\mathrm{a}}$. From Table \ref{tableMR}, we must
have
\[ \lambda = p + \frac{k}{2} p (p - 1), \]
for a certain $p \in \Z$. In this case $x_1$ is algebraic over $K
=\C(z)$ and $x_2 = x_1  \int x_1^{-2}$ is transcendental.

\begin{definition}
  Let $f (z)$ be a multivalued  function of the complex variable $z$,  and let $z_0 \in
  \mathbb{P}^1$.  We say that $e\in\C$ is the exponent of $f$ at $z_0$,  if in a neighbourhood of $z_0$, $f$ can be expressed into
  the  following form
  \[ f (z) = \zeta^e h (\zeta), \]
  where $\zeta$ is a local parameter around $z_0$,  $\zeta
  \mapsto h (\zeta)$ is holomorphic at $\zeta = 0$ and $h (0) \neq 0$.

  The principal part of $f$ at $z_0$ is denoted $f_{z_0}$, i.e, $f_{z_0} =
  \zeta^e h (0)$.

  We denote by $\mathcal{M}_{z_0}$ the monodromy operator around $z_0$.
\end{definition}

\begin{lemma}
\label{lem:ga}
  If  $G_1^{\circ} \simeq G_{\mathrm{a}}$ then,
  \begin{enumerate}
    \item Up to a complex multiplicative constant, the algebraic solution
    $x_1$ may be written in the form
\[ x_1 = z^a (z - 1)^b J (z) \mtext{where}a
    \in \left\{ \frac{k - 1}{2 k}, \frac{k + 1}{2 k} \right\}, \quad b \in \left\{ \frac{1}{4},
    \frac{3}{4} \right\},
\]
 and $J (z) \in \mathbb{R}[z]$ does not vanish at $z \in
    \{0, 1\}$.

    \item The function $\psi = \int x_1^{-2}$ has the exponent
    $1 - 2 b$ at $z = 1$ and $\mathcal{M}_1 (\psi) = - \psi$.
  \end{enumerate}
\end{lemma}

\begin{proof}
  (1)~For all $\sigma \in G_1$, $\sigma (x_1) = \chi (\sigma) x_1 + \mu
  (\sigma) x_2$ for certain $(\chi (\sigma) , \mu (\sigma)) \in\C^2$.
  But $\sigma (x_1)$ is still algebraic, hence $\mu (\sigma) = 0$, and $\sigma
  (x_1) = \chi (\sigma) x_1$. In particular $x_1$ is an eigenvector of the
  monodromy operators $\mathcal{M}_0$ and $\mathcal{M}_1$. For $|k| \geq
  3$, from equation \eqref{eqdeltaexp}, the differences of exponents at $z =
  0$ and $z = 1$ are not integers, hence we can deduce that $x_1$ is a
  principal branch of the Riemann scheme \eqref{eqRS} at $z = 0$ and at $z = 1$.
  Therefore, $x_1$ may be written in the form $x_1 = z^a (z - 1)^b J (z)$
  where $a$, (resp. $b$) is an exponent at $z = 0$ (resp. at $z = 1$), and $J (z)$
  is holomorphic on $\C$. Since $a$ and $b$ are rational numbers, $J
  (z) = x_1 / z^a (z - 1)^b$ is an algebraic function which is holomorphic on
  $\C$ hence, $J (z)$ is a polynomial. For $|k| = 1$, we have $\{ \frac{k -
  1}{2 k}, \frac{k + 1}{2 k} \}=\{0, 1\}$, therefore $x_1$ is regular at $z =
  0$. At $z = 1$ the difference of exponents is $\Delta_1 = 1 / 2$, so the
  previous arguments apply,  and point 1 is still true with $a \in \{0, 1\}$.
  Moreover, $J (0) \neq 0$ and $J (1) \neq 0$. Since the exponents are real, $J$ is a solution of a second order differential equation over $\R$. Thus, we can assume that
  $J\in\mathbb{R}[z]$.

  (2)~The function $x_1^{-2}$ has the exponent $- 2 b$ at $z = 1$. Thus,
  expanding it around $z=1$ and integrating, we obtain that $\psi$ has the
  exponent $1 - 2 b$ at $z = 1$. Therefore, $\mathcal{M}_1 (\psi) = \exp[ 2 \pi\rmi
  (1 - 2 b)] \psi = \exp[- 4 \pi  \rmi b] \psi = - \psi$,  because $b \in \{
  \frac{1}{4}, \frac{3}{4} \}$.
\end{proof}

Now, thanks to Remark~\ref{jordanredGa}, at  first we  have to test conditions
$(\alpha)$ for $\varphi$ and $\psi$. If we set
\[ \theta := z^a (z - 1)^b, \]
then, using  Lemma~\ref{lem:ga},  we have the explicit formulae
\[ \varphi = \int sx_1^2 = \int s \theta^2 J^2 = \frac{1}{2 k} \int z^{2 a -
   1} (z - 1)^{2 b - 1} J^2 (z) dz, \]
\[ \psi = \int x_1^{-2} = \int \frac{1}{\theta^2 J^2} . \]

\subsection{Algebraicity of $\psi$ and $\varphi$}

Since $G_1^{\circ} \simeq G_{\mathrm{a}}$, we know that $\psi$ is not algebraic and we
have the following.
\begin{lemma}
  \label{GaAlg} Let $|k| \geq 3$. If condition $(\alpha)$ holds
  then $\varphi$ is algebraic.
\end{lemma}

\begin{proof}
  Let $\tilde{L} =\C(z) [\theta^2]$ where $\theta = z^a (z - 1)^b$, and
\[ a
  \in \left\{ \frac{k - 1}{2 k}, \frac{k + 1}{2 k} \right\}, \mtext{and} b \in \left\{ \frac{1}{4},
  \frac{3}{4}\right\}.
\]
This is an algebraic extension of $K =\C(z)$ of
  degree
  \[ N =  \begin{cases}
       |k|  &\text{ when } k \in 2\mathbb{N},\\
       2| k| &\text{ when } k \not\in 2\mathbb{N}.
     \end{cases}
 \]
  Indeed, the minimal equation for $\theta^2$ is $(\theta^2)^N = z^{2 Na} (z -
  1)^{2 Nb} \in\C[z]$. Therefore, a basis of $\tilde{L} / K$ is $\{
  \theta^{-2}, 1, \theta^2, \cdots, (\theta^2)^{N - 2} \}$, and $N - 2
  \geq 2$ since $|k| \geq 3$. As $\varphi' \in \tilde{L}$, and
  $\psi' \in \tilde{L}$, from the Ostrowski-Kolchin theorem (see point 1 of Lemma \ref{ostrowski}),
  we deduce that condition $(\alpha)$ holds iff there exists $c \in
 \C$ such that, $\varphi + c \psi \in \tilde{L}$. But $\varphi + c
  \psi \in \tilde{L}$ iff there  exists a family $(f_{- 1}, \cdots, f_{N -
  2}) \in\C(z)^N$ such that
  \[ \varphi + c \psi = \sum_{i = - 1}^{N - 2} f_i ( \text{} \theta^2)^i . \]
  Differentiating the above equality, we obtain
  \begin{eqnarray*}
    \varphi' + c \psi' & = & \sum_{i = - 1}^{N - 2} \left(f_i ' + 2 i
    \frac{\theta'}{\theta} f_i\right) \theta^{2 i},\\
    \frac{1}{2 kz (z - 1)} \theta^2 J^2(z) + \frac{c}{\theta^2 J^2(z)} & = & \sum_{i
    = - 1}^{N - 2}\left(f_i ' + 2 i \frac{\theta'}{\theta} f_i\right) \theta^{2 i} .
  \end{eqnarray*}
  From this equation, we necessarily have
  \[
 \begin{split}
       \frac{c}{\theta^2 J^2}  = \left(f_{- 1} ' - 2 \frac{\theta'}{\theta} f_{-
       1}\right) \frac{1}{\text{} \theta^2} &\quad\Longleftrightarrow\quad   c \psi  =
       \frac{f_{- 1}}{\text{} \theta^2},  \\[0.5em]
       \frac{1}{2 kz (z - 1)} \theta^2 J^2  =  \left(f_1 ' + 2
       \frac{\theta'}{\theta} f_2\right) \theta^2 &\quad \Longleftrightarrow\quad \varphi  =
       f_2  \text{} \theta^2.
     \end{split}
\]
  The first equation implies that $c = 0$ because $\psi$ is not algebraic. The
  second equation implies that  $\varphi$ is algebraic. Moreover, $\varphi$ is
  algebraic iff there exists $f \in\C(z)$ such that
  \begin{equation}
    \frac{J^2(z)}{z (z - 1)} = f' + 2 \frac{\theta'}{\theta} f. \label{abeqphi}
  \end{equation}

\end{proof}

Since $\theta = z^a (z - 1)^b$, equation~\eqref{abeqphi} is equivalent to $J^2 = T (f)
:= z (z - 1) f' + 2 ((a + b) z - a) f$. Therefore we have the equivalence
\begin{equation}
  \varphi = \int sx_1^2 = f \theta^2 \quad \Longleftrightarrow\quad  J^2 = T (f) = z (z - 1)
  f' + 2 ((a + b) z - a) f. \label{eqphiT}
\end{equation}

\subsection{Algebraicity of $\varphi$ and condition $(\alpha)$}

At the end of the previous subsection we showed that $\varphi$ is algebraic iff
 the equation $J^2 = T (f)$ defined by~\eqref{eqphiT}, has a rational solution $f$.
The next Lemma gives an answer to this problem.
\begin{lemma}
  \label{GaPoly} Let $J \in \mathbb{R}[z]$ such that $J (0) J (1)
  \neq 0$. Then,
  \begin{enumerate}
     \item If $a \neq 1$, then the equation $J^2 = T (f)$ does not have rational
    solutions and $\varphi$ is not algebraic.

    \item If $a = 1$, and the equation $J^2 = T (f)$ has a solution $f \in
   \C(z)$, then  $f (z) = c (
   z^{-2} + 2 b z^{-1})+ g (z)$ where $c \neq 0$ is a constant, and
    $g (z)$ is a polynomial.

  \end{enumerate}
\end{lemma}

\begin{proof}
Let $f\in\C(z)$ be such that $J^2=T(f)$, in particular $T(f)$ is a polynomial.
We separate into three steps our further reasoning.

\noindent
\textit{First step.} We prove that $f$ has  only few poles,  precisely we claim that
\begin{enumerate}
 \item  if $ a\neq 1 $, then $ f \in \R[z]$;
\item  if $a=1$, then
 $ f (z) = c (
   z^{-2} + 2 b z^{-1})+ g (z)$  with $c\in\C$ and $\in\C[z]$.
\end{enumerate}
Indeed,   if $f$ has a pole of order $n$ at $t$, setting $f_t = c(z - t)^{-n}$,
  we have the following possibilities for the principal part of $T (f)$:
  \[
\begin{split}
       T (f)_t = \frac{- cnt (t - 1)}{(z - t)^{n + 1}}&\mtext{for} t\not\in\{0,1\}, \\
      T (f)_0 = \frac{c (n - 2 a)}{z^n}&\mtext{for} t=0,  \\
      T (f)_1 = \frac{c (2 b - n)}{(z - 1)^n}& \mtext{for} t=1.
     \end{split}
 \]
  If $t \not\in \{0, 1\}$, then $T (f)_t \neq 0$, so  $t$ is not a pole of $f$.  Similarly, since $2 b - n \neq 0$, $T(f)_1\neq 0$, and $t = 1$ cannot be a pole of $f$ .  Now, the formula $T (f)_0 = c (n - 2 a)z^{-n}$ is
  valid iff $n - 2 a \neq 0$. But $n - 2 a = n - 1 \pm \frac{1}{k}$, thus
  \[ (n - 2 a = 0 )\quad \Longleftrightarrow\quad \left(n = 1 \mp \frac{1}{k} \right) \quad\Longleftrightarrow \quad
     \left(  n = 2  \mtext{and}
       a = 1 \mtext{and}
       k = \mp 1\right).
\]
  Therefore, if $a \neq 1$, then $f$ does not have pole at  $z = 0$, and $f$ must be a polynomial.
   Now  if  $f$ has a pole at $z=0$,   according  to the previous equivalence, we must have $a=1$ and $n=2$.
But,  if $a = 1$, then $T (f) = z (z - 1) f' + 2 ((1 + b) z - 1) f$, and we have
  \[
T \left( \frac{1}{z^2}\right) = \frac{2 b}{z} \mtext{ and } T \left( \frac{1}{z}\right) =
     -\frac{1}{z} + 2 b + 1. \]
  If $f$ is a solution which is not a polynomial, it must have a pole of order
  two at zero, and, for the compensation, we must have $f (z) = c ( z^{-2} +
  2 bz^{-1}) + g (z)$, where $c \in\C$ and $g (z) \in
 \C[z]$.

\noindent
 \textit{Second step.} We now treat the particular case $a=0$.   If $f$ is a
  rational solution of the equation $J^2 = T (f)$, then, by the first step, $f$ is a polynomial.
  Evaluating this equation at $z = 0$ we get
  \[ J^2 (0) = - 2 af (0) . \]
  Therefore, if $a = 0$, then  $J (0) = 0$,  but it is not true. Thus, in this case,
   the equation does not have rational solutions.

\noindent
 \textit{Third step.} Under the assumption that    $a \neq 0$  we claim that the equation $T(f)=J^2$  does not have  polynomial solutions.  Since  $\theta = z^a (z - 1)^b $, equivalence
  \eqref{eqphiT} can  be written in the following form
  \[ z^{2 a - 1} (z - 1)^{2 b - 1} J^2 (z) = \frac{d}{dz} (f (z) z^{2 a} (z -
     1)^{2 b}) ,
 \]
where
\[
 a \in \left\{ \frac{k - 1}{2 k},
  \frac{k + 1}{2 k} \right\}, \quad b \in\left \{ \frac{1}{4}, \frac{3}{4} \right\}.
\]
Hence, since $a \neq 0$, we have $2 a = 1 \pm \frac{1}{k} > 0$, and moreover, $2 b \geq \frac{1}{2}$.  Therefore integrating between 0 and 1 we
  get
  \[
  f (z) z^{2 a} (z - 1)^{2 b}\Big|_0^1 = 0 = \int z^{2 a - 1} (z - 1)^{2 b -
     1} J^2 (z) dz > 0,
 \]
  since the integrand is positive. The above contradiction proves the claim.  As a conclusion, if $a \neq 1$, the equation $J^2 = T (f)$ does not have
  rational solution. This proves Point 1. When $a = 1$, and $J^2 = T (f)$
  possesses a rational solution, the latter cannot be a polynomial. and by the
  first step, point 2 follows.
\end{proof}

In the case $a = 1$, which  happens only for  $k = \pm 1$,  $\varphi$ can be algebraic,  so condition $(\alpha)$ can be
satisfied. For example, computations
with Riemann schemes show that for row 4 in Table \ref{tableMR}, when $(k,
\lambda) = (- 1, - 2)$, we have $x_1 = z (z - 1)^{3 / 4}$, and
\[ \varphi = \int z (1 - z)^{1 / 2} dz = \frac{6 z^2 - 2 z - 4}{15}  \sqrt{1 -
   z} \in L =\C(z) [x_1] =\C(z) [(1 - z)^{1 / 4}], \]
is therefore algebraic and condition $(\alpha)$ is satisfied. Nevertheless,  for those cases we have the following.
\begin{lemma}
  \label{alpha}Let us assume that $a = 1$. If condition $(\alpha)$ holds,
  then condition $(\beta)$ is not satisfied.
\end{lemma}
\begin{proof}
  By Definition \ref{abcondition} and Theorem \ref{abcritere},  we have
  to check if the  condition
  \[ Q = \varphi \psi - 2   \int \varphi \psi'  \in L[\psi] ,\]
is satisfied.
  By assumption,   $\varphi + c \psi\in L$,  for a certain
 $c\in\C$. Thus, we have
  \[ \left(Q \in L [\psi]\right) \quad \Longleftrightarrow\quad \left( I (z) = \int \psi' \cdot \varphi \in L
     [\psi] \right),\]
  where
  \[ I (z) = \int \psi' \cdot \varphi = \int \frac{f \theta^2}{J^2 \theta^2} =
     \int \frac{f (z)}{J^2 (z)} dz. \]
 As  $a = 1$, by point 2 of  Lemma \ref{GaPoly},  $f (z) = c (z^{-2} +
  2 bz^{-1})+ g (z)$ with $c \neq 0$. Therefore $I (z)$ may be expressed
  by a formula of the form
  \[
I (z) = \gamma_0 \mathrm{Log} (z) + \sum \gamma_i \mathrm{Log} (z - z_i) + h
     (z),
 \]
  where, $h (z) \in\C(z)$, $\gamma_0 = {- 2 bc}/{J^2 (0)} \neq
  0$,  $\gamma_i \in\C$, and  $z_i$ are roots of $J (z)$. In
  particular $z_i \not\in \{0, 1\}$. Hence, $I (z) \in L [\psi]$ if and only if
  $I (z)$ and $\psi (z)$ are algebraically dependent. But, by the the Ostrowski-Kolchin
  theorem, this happen if and only if we have a non trivial linear relation
  with complex coefficients
  \[ \mu I (z) + \nu \psi (z) = \omega (z) \in L. \]
  However, $\mathcal{M}_1 (I (z)) = I (z)$ and, from Lemma \ref{lem:ga},
  $\mathcal{M}_1 (\psi) = \exp[-\pi\rmi] \psi = - \psi$. Applying the
  monodromy operator to the previous equation yields
  \[ \mu I (z) - \nu \psi (z) =\mathcal{M}_1 (\omega (z)) . \]
  So, $2 \mu I (z) = \omega (z) +\mathcal{M}_1 (\omega (z))$ is algebraic. As
  $I (z)$ is not algebraic, because $\gamma_0 \neq 0$, we deduce that $\mu = \nu
  = 0$ and \ condition $(\beta)$ is not satisfied.
\end{proof}

\subsection{Proof of Theorem \ref{jordanGa}}

  By Remark \ref{jordanredGa}, it is enough to show that  $G_2^{\circ}$
  is not Abelian. Since here $G_1^{\circ} \simeq G_{\mathrm{a}}$, from Theorem
  \ref{abcritere}, it remains to show that conditions $(\alpha)$ and
  $(\beta)$ are not simultaneously satisfied.

  From Remark \ref{jordanredGa} again, the pair $(k, \lambda)$ must belong
  to rows 2, 3, or 4 of Table \ref{tableMR}. In particular, either $|k| \geq
  3$, or  $k = \pm 1$.
  \begin{itemize}
    \item For $|k| \geq 3$, condition $(\alpha)$ is not satisfied.
    Indeed, from Lemma \ref{GaAlg}, condition $(\alpha)$ implies that $\varphi$ is
    algebraic and, from point 1 of Lemma \ref{GaPoly}, we know that in this
    case $\varphi$ is not algebraic.

    \item For  $k = \pm 1$,  condition $(\alpha)$ may be satisfied but if this happens,
    by Lemma \ref{alpha}, condition $(\beta)$ is not satisfied.
  \end{itemize}
The above finishes the proof.\qed

\section{Elimination of the Jordan blocks with $ G_1^{\circ} \simeq
\{\id\}$}

Our main result in this section is the following.
\begin{theorem}
  \label{abfiniteG1}Assume that $V'' (c)$ has a Jordan block of
  size $d \geq 3$, and $G_1$ is a finite subgroup of $\mathrm{SL}(2, \C)$.  Then $G (\mathrm{VE}_z)^{\circ}$ is not Abelian. This
   eliminates  the rows with numbers from 5 to 21 in Table \ref{tableMR}.
\end{theorem}

\begin{remark}
  \label{jordanredId}
Let $B (\lambda, d)$ be a Jordan block of
  $V'' (c)$ with size $d \geq 3$ and eigenvalue $\lambda$. Since $G_1^{}$
  is finite and correspond to the VE
  \[ \frac{d^2 \eta}{dt^2} = - \lambda \varphi^{k - 2} (t) \eta, \]
  we deduce from Theorem \ref{MR}, that necessarily, the pair $(k, \lambda)$
  must belong to rows 5 to 21 of Table~\ref{tableMR}. Now, passing to the VE
  in the $z$ variable, we know that the system
  \[ \begin{bmatrix}
       x''\\
       y''\\
       u''
     \end{bmatrix}= \begin{bmatrix}
       r & 0 & 0\\
       s & r & 0\\
       0 & s & r
     \end{bmatrix} \begin{bmatrix}
       x\\
       y\\
       u
     \end{bmatrix} , \]
  with Galois
  group $G_3$ is a subsystem of VE (\ref{eqVEredz}). From Proposition \ref{VEtVEz}, $G
  (\mathrm{VE}_t)^{\circ} \simeq G (\mathrm{VE}_z)^{\circ}$ and $G_3^{\circ}$ is a
  quotient of $G (\mathrm{VE}_z)^{\circ}$. Therefore it is enough to prove
  that $G_3^{\circ}$ is not Abelian.
\end{remark}

 Recall from Lemma
\ref{abkovacic} that if $G_1$ is finite, then it is  one of the following types
\begin{enumerate}
  \item Dihedral type: $G_1$ is conjugated to a finite subgroup of
\begin{equation*}
 D^{\dag} = \defset{
\begin{bmatrix}
    \lambda & 0\\
       0 & 1 / \lambda
\end{bmatrix}
}{\lambda\in\C^\star} \cup
 \defset{
\begin{bmatrix}
   0 & \lambda\\
       - 1 / \lambda & 0
\end{bmatrix}
}{ \lambda\in\C^\star}
\end{equation*}

  \item Tetrahedral type: $G_1 /\{\pm \id\} \simeq \mathfrak{A}_4$

  \item Octahedral type: $G_1 /\{\pm \id\} \simeq \mathfrak{S}_4$

  \item Icosahedral type: $G_1 /\{\pm \id\} \simeq \mathfrak{A}_5$
\end{enumerate}
From Theorem \ref{abcritere}, we know that if $G_1$ is finite, $G_2^{\circ}$
is Abelian, where $G_2$ is the Galois group of the two first equations of the
above system. This why we have to consider  Jordan blocks of size $d \geq
3$, in order to find obstructions to the integrability.
At first, we  build   some theoretical results in the spirit of
Section 2.

\subsection{Theory for Jordan blocks of size three}

Now we assume that the size of the Jordan block is three. With the notations
of Section 2, the subsystem of the variational equations corresponding to the block,  can be written in the
following two equivalent forms
\begin{equation}
  \left\{\begin{matrix}
    x'' =& rx& &\\
    y'' =& ry &+& sx\\
    u'' =& ru &+& sy
  \end{matrix} \right. \quad \Longleftrightarrow\quad \begin{bmatrix}
    x'\\
    x''\\
    y'\\
    y''\\
    u'\\
    u''
  \end{bmatrix}= \begin{bmatrix}
    R & 0 & 0\\
    S & R & 0\\
    0 & S & R
  \end{bmatrix}\begin{bmatrix}
    x\\
    x'\\
    y\\
    y'\\
    u\\
    u'
  \end{bmatrix},\label{abvar3}
\end{equation}
where $R$ and $S$ are $2\times 2$ matrices given by \eqref{eq:RS}.

Let us fix more  notations.
  \begin{itemize}
    \item $F_1 / K$ is the Picard-Vessiot extension associated to the equation
    $L_2 (x) = x'' - rx = 0$. Its Galois group is still denoted by $G_1$.

    \item $F_2 / K$ is the Picard-Vessiot extension associated to the first
    two equations of \eqref{abvar3}. Its Galois group, is still denoted by
    $G_2$.

    \item $F_3 / K$ is the Picard-Vessiot extension over $K$ associated to
    \eqref{abvar3}. Its Galois group is denoted by $G_3$.
  \end{itemize}

\begin{remark}
  We have the following inclusions of differential fields
  \[ K \subset F_1 \subset F_2 \subset F_3 . \]
  All the results of Section 2 can be applied to the extension $F_2 / K$. In
  particular, since $G_1$ is finite, from Theorem \ref{abcritere},
  $G_2^{\circ}$ is Abelian. Therefore, $G_2^{\circ}$ is an Abelian quotient of
  $G_3^{\circ}$.
\end{remark}

We fix a basis $\{x_1, x_2 \}$ of the solution space $V$ of $L_2 = 0$.  Let
$(y_1, y_2, u_1, u_2)$ be an element of $F_3^4$ such that
\[
\left\{ \begin{matrix}
     y_1'' = ry_1 + sx_1\\
     y_2'' = ry_2 + sx_2
   \end{matrix} \right. \mtext{ and } \left\{\begin{matrix}
     u_1'' = ru_1 + sy_1\\
     u_2'' = ru_2 + sy_2
   \end{matrix} \right. . \]
Then, we set
\[ X = \begin{bmatrix}
     x_1 & x_2\\
     x_1' & x_2'
   \end{bmatrix}, \quad  Y = \begin{bmatrix}
     y_1 & y_2\\
     y_1' & y_2'
   \end{bmatrix}, \quad U = \begin{bmatrix}
     u_1 & u_2\\
     u_1' & u_2'
   \end{bmatrix}, \quad \Xi_3 = \begin{bmatrix}
     X & 0 & 0\\
     Y & X & 0\\
     U & Y & X
   \end{bmatrix}. \]
Similarly as in Section 2, $\Xi_3$ is a fundamental matrix of solutions of
\eqref{abvar3}.

For all $\sigma \in G_3$, the equation $\sigma (\Xi_3) = \Xi_3 M (\sigma)$,
forces $\sigma$ to be represented by a $6 \times 6$ matrix $M (\sigma)$ of
the form
\[ M (\sigma) = \begin{bmatrix}
     A (\sigma) & 0 & 0\\
     B (\sigma) & A (\sigma) & 0\\
     C (\sigma) & B (\sigma) & A (\sigma)
   \end{bmatrix}. \]
\begin{proposition}
  \label{abcriteredim3}Assume that  $G_1$ is finite. Then $G_3^{\circ}$ is
  Abelian iff there exists a basis $\{x_1 , x_2 \}$ of $V = \mathrm{Sol} (L_2)$
  such that  one of the following condition is satisfied
  \begin{itemize}
    \item $\varphi_1 = \int sx_1^2 \in F_1$ and $\int \varphi'_1 \psi_1 \in
    F_1$ where $\psi_1 = \int x_1^{-2}$.

    \item $\int sx_1^2 \in F_1$ and $\int sx_2^2 \in F_1$.
  \end{itemize}
  If $G_3^{\circ}$ is Abelian, then there exists at least one non-zero $x \in
  V = \mathrm{Sol} (L_2)$ such that $\int sx^2 \in F_1$.
\end{proposition}
\begin{proof}
We consider  $G_3^{\circ}$  as a subgroup of $\mathrm{SL}(6 ,\C)$. The   elements of $G_3^{\circ}$
  are matrices the form
  \[ P (B , C) := \begin{bmatrix}
       \id & 0 & 0\\
       B & \id & 0\\
       C & B & \id
     \end{bmatrix}. \]
  The product and the commutators of two such matrices are given by
  \begin{gather*}
 P(B_1 , C_1) P (B_2 , C_2) = P(B_1 + B_2 , C_1 + C_2 + B_1 B_2),\\
 [P(B_1 , C_1) P (B_2 , C_2)]=P(0,[B_1,B_2]).
\end{gather*}
Set
  \[
\mathbb{B}:=\defset{B \in \M (2, \C)}{\exists C \in \M (2,\C)
     \mtext{and}P (B , C) \in G_3^{\circ} }.
\]
  Then, thanks to the above formulae, $G_3^{\circ}$ is Abelian iff any two
  matrices belonging to $\mathbb{B}$ commute. This is the case iff, up to
  conjugation,  $\mathbb{B}$ is  contained either in the set of upper
  triangular matrices with diagonal of the form $a \id$, or,
  $\mathbb{B}$ is contained in the set of diagonal matrices. For any of
  this two cases, thanks to a conjugation formula similar to~\eqref{eqab7}, we can find a basis  $\{x_1 , x_2 \}$  of $V$ such that the
  representation of the elements of $\mathbb{B}$ in this basis are either
  upper triangular or diagonal.

  From point 4 of Proposition \ref{abG},  we have $G_2^{\circ} = H =
  \mathrm{Gal}_{\partial} (F_2 / F_1)$. Let
  \[ \pi_2 : G_3^{\circ} \rightarrow G_2^{\circ},\quad  P (B , C) \mapsto
     \begin{bmatrix}
       \id & 0\\
       B & \id
     \end{bmatrix}=N(B), \]
be the projection.
  With the notations of Proposition \ref{abcritereH}, the two above conditions
  for $\mathbb{B}$ are respectively equivalent to $\pi_2 (G_3^{\circ})
  = G_2^{\circ} \subset H_{\mathrm{a}}$,  and  $\pi_2 (G_3^{\circ}) = G_2^{\circ} \subset
  H_{\mathrm{m}}$.

  But,  from Proposition \ref{abcritereW}, we have
  \[ \left\{ \begin{matrix}
       G_2^{\circ} \subset H_{\mathrm{a}} & \ \Longleftrightarrow \  & W_1 \text{ and } W_3 \text{
       hold},\\
       G_2^{\circ} \subset H_{\mathrm{m}} & \  \Longleftrightarrow\  & W_1 \text{ and } W_2 \text{
       hold}.
     \end{matrix} \right.  \]
  Now, from Definition \ref{defW} and Lemma \ref{abcomputation},  condition $W_1$ holds
  iff $W (x_1 , y_1) =\varphi_1= \int sx_1^2 \in F_1$, and the same result holds for
  condition $W_2$. From the same definition and lemma, condition $W_3$ holds iff \[ Q = \varphi_1
  \psi_1 - 2 \int \varphi_1 \psi_1' =  - \varphi_1
  \psi_1+2 \int \varphi_1' \psi_1 \in F_1.\]
 But $\psi_1 = x_2 / x_1 \in F_1$,  and $\varphi_1$ also
  belongs to $F_1$ if $W_1$ is assumed to be satisfied. Therefore,  $W_3$ holds
  iff $\int \varphi'_1 \psi_1 \in F_1$
\end{proof}
\begin{proposition}
  \label{abfinite} Let $V = \mathrm{Sol} (L_2)$.
Assume that $G_1$ is finite. Then we have
  the following properties
  \begin{enumerate}
    \item Let $x_1$ be a non-zero element of $V$. If $\int sx_1^2 \in F_1$, then for all $\sigma \in G_1$, $\int s
    \sigma (x_1)^2 \in F_1$.

    \item For all $x \in V$, $\int sx^2 \in F_1$, iff there exists a basis
    $\{x_1, x_2 \}$ of $V$ such that
    \[ \int sx_1^2 \in F_1 \text{ and } \int sx_1 x_2 \in F_1 \text{ and }
       \int sx_2^2 \in F_1 . \]
    \item Assume that $\int sx_1^2 \in F_1$, and $\int sx^2 \not\in F_1$,  then
    $G_1$ is of dihedral type.

  \end{enumerate}
\end{proposition}

\begin{proof}
  $(1)$ Let $x_1$ be any non zero element of $V$.  Since $\int sx_1^2 = W (x_1
  , y_1) \in F_2$, for all $\sigma \in G_2$ we have
 \[ \sigma ( \int sx_1^2) = \sigma
  (W (x_1 , y_1)) = W (\sigma (x_1) , \sigma (y_1)) = \int s \sigma (x_1)^2.
\]
  Therefore, if $\int sx_1^2 \in F_1$, then $\sigma ( \int sx_1^2) = \int s
  \sigma (x_1)^2 \in F_1$. Since $F_1 / K$ is a Picard-Vessiot extension
  contained in $F_2 / K$,  the restriction morphism $\mathrm{Res} : G_2
  \rightarrow G_1$ is surjective, therefore the integrals $\int s \sigma
  (x_1)^2 \in F_1$ for all $\sigma \in G_1$.

  $(2)$ Assume that for all $x \in V$, $\int sx^2 \in F_1$, and let $\{x_1, x_2
  \}$ be a basis of $V$.  Then the three particular integrals
  \[ \int s (x_1 + x_2)^2 \text{ and } \int sx_1^2 \text{ and } \int sx_2^2,
  \]
  belong to $F_1$. Taking the difference of those integrals we deduce that $\int
  sx_1 x_2 \in F_1$. Conversely, each $x \in V$ can be written in the form
  $x = \lambda x_1 + \mu x_2$. Therefore,
  \[ \int sx^2 = \lambda^2  \int sx_1^2 + 2 \lambda \mu \int sx_1 x_2 + \mu^2
     \int sx_2^2 \in F_1 . \]
  (3) For the action of $G_1$ on $\mathbb{P}(V) \simeq \mathbb{P}^1$, when
  we look at the orbit $\Omega$ of $[x_1]$, three cases may a priory happen:
  \begin{itemize}
    \item[a)] $\mathrm{Card} (\Omega) = 1$.

    \item[b)] $\mathrm{Card} (\Omega) = 2$.

    \item[c)] $\mathrm{Card} (\Omega) \geq 3$.
  \end{itemize}
  Let us  first prove that with the assumption of point 3, case c) cannot
  happen. Indeed case c) implies that there exists $x_2 = \sigma_1 (x_1)$ which
  is not collinear to $x_1$, and also there exists $x_3 = \sigma_2 (x_1) =
  \lambda x_1 + \mu x_2$ with $\lambda \mu \neq 0$. From point 1, this
  implies that the three integrals
  \[ sx_1^2 \mtext{ and } \int sx_2^2 \mtext{ and } \int s (\lambda x_1 + \mu
     x_2)^2, \]
  belong to $F_1$. So, $\int sx_1 x_2$  belongs to $F_1$. Thus, from point 2,
  for all $x \in V$, $\int sx^2 \in F_1$ which is not true. There remains to
  show that in cases a) and b), $G_1$ is of dihedral type.

  In case b), let $\Omega =\{[x_1], [x_2]\}$. This means that $\{x_1, x_2 \}$
  is a basis of $V$. Moreover, any conjugate of $x_1$ or $x_2$ is either
  collinear to $x_1$, or to $x_2$. Hence, in the basis $\{x_1, x_2 \}$, the
  representation of $G_1$ is of dihedral type.

  In Cases a), since $\Omega =\{[x_1]\}$, $x_1$ is an eigenvector of any
  $\sigma \in G_1$. We find a second common eigenvector for any $\sigma \in
  G_1$,  using the following classical averaging argument coming from
  representation theory. Let $\langle \cdot, \cdot\rangle$ be an arbitrary Hermitian product on $V \simeq
 \C^2$ for which $x_1$ is not an isotropic vector (i.e. $\langle x_1, x_1
  \rangle \neq 0$). Consider the average
  \[ (X, Y) = \sum_{\sigma \in G_1} \langle\sigma (X), \sigma (Y) \rangle . \]
  The pairing $(\cdot,\cdot)$ is a new Hermitian product on $V$ for which $G_1$ is
  unitarian. Therefore the orthogonal of the line $\C x_1$ is another
  line of the form $\C x_2$ which is also globally $G_1$-invariant.
  Therefore $G_1$ is diagonalizable in the basis $\{x_1, x_2 \}$. This
  proves that $G_1$ is of dihedral type.

\end{proof}
\begin{proposition}
  \label{abeuler}Assume that $G_1$ is finite, $K =\C(z)$, $s
  = \frac{1}{2 kz (z - 1)}$, and consider the following properties
  \begin{enumerate}
    \item For all $x \in V$, $\int sx^2 \in F_1$.

    \item $F_2 = F_1$ and $G_2 \simeq G_1$.

    \item There exists $M \in \mathrm{GL}(2,K)$ such that $S = M' + [M , R]$,
    where $R$ and $S$ are given by \eqref{eq:RS}.

    \item There exists a non-zero rational solution $v \in\C(z)$ to the
    equation
    \[ L_4 (v) = [z (z - 1) L_2^{\circledS 2} (v)]' = 0, \]
     where $L_2^{\circledS 2}$ denotes the second symmetric power of
    $L_2$.
  \end{enumerate}
  Then we have $(1) \Leftrightarrow (2) \Leftrightarrow (3)\Rightarrow
  (4)$.
\end{proposition}

\begin{proof}
  $(1) \Leftrightarrow (2)$ Let $\{x_1 , x_2 \}$ be a basis of $V$. From point
  $2$ of Proposition \ref{abfinite}, property 1 is equivalent to
  \[ W (x_1 , y_1) \in F_1 \mtext{ and } W (x_1 , y_2) \in F_1 \mtext{ and } W
     (x_2 , y_2) \in F_1 . \]
  By Proposition \ref{abcritereW} and Corollary \ref{abtrace}, these three
  Wronskians are fixed by the elements
\[
\sigma = \begin{bmatrix}
    \id & 0\\
    B (\sigma) & \id
  \end{bmatrix}\in \mathrm{Gal}_{\partial} (F_2 / F_1)
\]
iff $B (\sigma) =
  0$. So,   property 1 is equivalent to  $\mathrm{Gal}_{\partial} (F_2 / F_1)
  =\{\id\}$, that is to property 2.

  $(2) \Leftrightarrow (3)$ From the exact sequence
  \[
\begin{CD}
       \{\id\} @>>> \mathrm{Gal}_{\partial} (F_2 / F_1) @>>> G_2 @>>> G_1 @>>> \{\id\}\\
     @. @. \begin{bmatrix}
         A (\sigma) & 0\\
         B (\sigma) & A (\sigma)
       \end{bmatrix}@>>> A (\sigma)  @.
\end{CD}
\]
  we have
\[ \left( \mathrm{Gal}_{\partial} (F_2 / F_1) =\{\id\} \right) \  \Longleftrightarrow \
 \left(  \forall \sigma \in G_2, B (\sigma) = 0 \right) .
\]
But the general formulae for the
  action of $G_2$ are $\sigma (X) = XA (\sigma)$ and $\sigma (Y) = YA (\sigma)
  + XB (\sigma)$. This implies that
  \[ \sigma (YX^{- 1}) = YX^{- 1} + XB (\sigma) A^{- 1} (\sigma) X^{- 1} . \]
  So, $\sigma (YX^{- 1}) = YX^{- 1}$ iff $B (\sigma) = 0$. Therefore
  $\mathrm{Gal}_{\partial} (F_2 / F_1) =\{\id\}$ iff $ YX^{- 1}
  \in \mathrm{GL}(2 ,K)$.

  Now we are looking for the differential equation satisfied by $M = YX^{-
  1}$. From
  \[ \Xi_2= \begin{bmatrix}
       X & 0\\
       MX & X
     \end{bmatrix}, \quad \Xi'_2 = \begin{bmatrix}
       X' & 0\\
       M' X + MX' & X'
     \end{bmatrix}= \begin{bmatrix}
       R & 0\\
       S & R
     \end{bmatrix}\Xi_2, \]
  we obtain
  \begin{eqnarray*}
    X' & = & RX,\\
    M' X + MX' & = & SX + RMX,\\
    M' X + MRX & = & SX + RMX,\\
    M' + MR - RM & = & S.
  \end{eqnarray*}
  This proves $(2) \Leftrightarrow (3)$.

  $(3) \Rightarrow (4)$. We write $M =
 \bigl[ \begin{smallmatrix}
    u & v\\
    f & g
  \end{smallmatrix}\bigr]$, and we insert this expression into the above
  differential equation. This gives a system of four equations. By expressing
  $f$ and $g$ in terms of $u$ and $v$, the original equation is equivalent to
  the system
  \[
\left\{
\begin{split}
       f & =  u' + rv,\\
       g & =  u + v',\\
       u' & =  - v'' / 2,\\
       s & =  u'' + r' v + 2 rv'.
     \end{split}\right.
 \]
  From the above, $v$ satisfies $L_3 (v) : = v''' - 4 rv' - 2 r' v = - 2 s$.
  But $L_3 (v) = L_2^{\circledS 2} (v)$ is the second symmetric power of $L_2
  (v) = v'' - rv$. Now, for $K =\C(z)$ and $s = \frac{1}{2 kz (z -
  1)}$,
  \[ L_3 (v) = - 2 s = \frac{- 1}{kz (z - 1)} \ \Longrightarrow \  L_4 (v) : = [z (z -
     1) L_2^{\circledS 2} (v)]' = 0. \]
  Hence, if $M \in \mathrm{GL}(2 ,\C(z))$, then $v \in\C(z)$,
  and this implies that the equation $L_4 (v) = 0$ has a non-zero rational solution.
\end{proof}
Surprisingly, the differential equation $S = M' + [M , R]$ has the form  the classical Euler
equation for the angular momentum of a  rigid body, see \cite{arnold89}
pp.142-143.

\subsection{Type of $G_1$ when it is finite}

In order to apply the previous theory, we need to compute $G_1 =\mathcal{G}(k,
\lambda)$ when $G_1^{\circ} =\{\id\}$ in Table \ref{tableMR}. Table~\ref{tablefiniteG} below gives this information.

\begin{table}[th]
 \caption{\label{tablefiniteG} Type of $\mathcal{G}(k,\lambda)$ with $\mathcal{G}(k,\lambda)^\circ=\{\id\}$.}
  \begin{tabular}{|c|c|c|c|c|}
\hline
row & k & $\lambda$ & Exponents of $L_2$ at $\{0, \infty\}$ &  $\mathcal{G}(k , \lambda)$ \\
    \hline
    5 & 1 & $0$ & $(0, 1), (- 1 / 4, - 3 / 4)$ & cyclic-dihedral\\
    \hline
    6 & -1 & $1$ & $(0, 1), (- 1 / 4, - 3 / 4)$ & cyclic-dihedral\\
    \hline
    7 & $|k| \geq 3$ & $\dfrac{1}{2}  \left( \dfrac{k - 1}{k} + p (p + 1)
    k \right)$ & $\varepsilon_0$, $\left( \dfrac{1}{4} (2 p - 1), - \dfrac{1}{4} (2
    p + 3)\right)$ & dihedral\\
    \hline
    8 & 3 & $\dfrac{- 1}{24} + \dfrac{1}{6} (1 + 3 p)^2$ & $\left( \dfrac{1}{3},
    \dfrac{2}{3}\right), \left(- \dfrac{2}{3} - \dfrac{1}{2} p, - \dfrac{1}{3} + \dfrac{1}{2}
    p\right)$ & tetrahedral\\
    \hline
    9 & 3 & $\dfrac{- 1}{24} + \dfrac{3}{32} (1 + 4 p)^2$ & $\left( \dfrac{1}{3},
    \dfrac{2}{3}\right),\left(- \dfrac{5}{8} - \dfrac{1}{2} p, - \dfrac{3}{8} + \dfrac{1}{2}
    p\right)$ & octahedral\\
    \hline
    10 &3  & $\dfrac{- 1}{24} + \dfrac{3}{50} (1 + 5 p)^2$ & $\left( \dfrac{1}{3},
    \dfrac{2}{3}\right), \left(- \dfrac{3}{5} - \dfrac{1}{2} p, - \dfrac{2}{5} + \dfrac{1}{2}
    p\right)$ & icosahedral\\
    \hline
    11 & 3 &  $\dfrac{- 1}{24} + \dfrac{3}{50} (2 + 5 p)^2$ & $\left( \dfrac{1}{3},
    \dfrac{2}{3}\right), \left(- \dfrac{7}{10} - \dfrac{1}{2} p, - \dfrac{3}{10} +
    \dfrac{1}{2} p\right)$ & icosahedral\\
    \hline
    12 & -3 & $\dfrac{25}{24} - \dfrac{1}{6} (1 + 3 p)^2$ & $\left( \dfrac{1}{3},
    \dfrac{2}{3}\right), \left(- \dfrac{2}{3} - \dfrac{1}{2} p, - \dfrac{1}{3} + \dfrac{1}{2}
    p\right)$ & tetrahedral\\
    \hline
    13 & -3 &  $\dfrac{25}{24} - \dfrac{3}{32} (1 + 4 p)^2$ & $\left( \dfrac{1}{3},
    \dfrac{2}{3}\right), \left(- \dfrac{5}{8} - \dfrac{1}{2} p, - \dfrac{3}{8} + \dfrac{1}{2}
    p\right)$ & octahedral\\
    \hline
    14 & -3 & $\dfrac{25}{24} - \dfrac{3}{50} (1 + 5 p)^2$ & $\left( \dfrac{1}{3},
    \dfrac{2}{3}\right), \left(- \dfrac{3}{5} - \dfrac{1}{2} p, - \dfrac{2}{5} + \dfrac{1}{2}
    p\right)$ & icosahedral\\
    \hline
    15 & -3 & $\dfrac{25}{24} - \dfrac{3}{50} (2 + 5 p)^2$ & $\left( \dfrac{1}{3},
    \dfrac{2}{3}\right), \left(- \dfrac{7}{10} - \dfrac{1}{2} p, - \dfrac{3}{10} +
    \dfrac{1}{2} p\right)$ & icosahedral\\
    \hline
    16 & 4 & $\dfrac{- 1}{8} + \dfrac{2}{9} (1 + 3 p)^2$ & $\left( \dfrac{3}{8},
    \dfrac{5}{8}\right), \left(- \dfrac{2}{3} - \dfrac{1}{2} p, - \dfrac{1}{3} + \dfrac{1}{2}
    p\right)$ & octahedral\\
    \hline
    17 & -4 & $\dfrac{9}{8} - \dfrac{2}{9} (1 + 3 p)^2$ & $\left( \dfrac{3}{8},
    \dfrac{5}{8}\right), \left(- \dfrac{2}{3} - \dfrac{1}{2} p, - \dfrac{1}{3} + \dfrac{1}{2}
    p\right)$ & octahedral\\
    \hline
    18 & 5 & $\dfrac{- 9}{40} + \dfrac{5}{18} (1 + 3 p)^2$ & $\left( \dfrac{2}{5},
    \dfrac{3}{5}\right), \left(- \dfrac{2}{3} - \dfrac{1}{2} p, - \dfrac{1}{3} + \dfrac{1}{2}
    p\right)$ & icosahedral\\
    \hline
    19 & 5 &  $\dfrac{- 9}{40} + \dfrac{1}{10} (2 + 5 p)^2$ & ($\dfrac{2}{5},
    \dfrac{3}{5}) , \left(- \dfrac{7}{10} - \dfrac{1}{2} p, - \dfrac{3}{10} +
    \dfrac{1}{2} p\right)$ & icosahedral\\
    \hline
    20 & -5 & $\dfrac{49}{40} - \dfrac{5}{18} (1 + 3 p)^2$ & $\left( \dfrac{2}{5},
    \dfrac{3}{5}\right), \left(- \dfrac{2}{3} - \dfrac{1}{2} p, - \dfrac{1}{3} + \dfrac{1}{2}
    p\right)$ & icosahedral\\
    \hline
    21 & -5 & $\dfrac{49}{40} - \dfrac{1}{10} (2 + 5 p)^2$ & $\left(\dfrac{2}{5},
    \dfrac{3}{5}\right) ,\left(- \dfrac{7}{10} - \dfrac{1}{2} p, - \dfrac{3}{10} +
    \dfrac{1}{2} p\right)$ & icosahedral\\
    \hline
  \end{tabular}
\end{table}

To determine the last column of Table \ref{tablefiniteG}, we used the
following facts.
\begin{itemize}
  \item The exponents of $L_2$ at $\{0, 1, \infty\}$ are
  \[ \varepsilon_0 = \left\{ \frac{k - 1}{2 k}, \frac{k + 1}{2 k} \right\}, \quad
     \varepsilon_1 = \left\{ \frac{1}{4}, \frac{3}{4} \right\}, \quad
     \varepsilon_{\infty} = \left\{ \frac{\tau - 1}{2}, - \frac{\tau + 1}{2}
     \right\}, \]
  with \[ \tau = \frac{1}{2 k} \sqrt{(k - 2)^2 + 8 k \lambda}.\]

  \item Therefore, the reduced differences of exponents are
  \[ \Delta_0 = \Bigl| \frac{1}{k} \Bigr|, \quad \Delta_1 = \frac{1}{2}, \quad \Delta_{\infty} = |
     \tau | \mod \Z.
\]
\end{itemize}
Thanks to the Schwarz Table, see p.128 in \cite{poole}, we can
compute   $G_1/\{\pm\id\} =\mathcal{G}(k, \lambda)/\{\pm\id\} $ which is the image of $G_1$ in $\mathrm{PSL}(2,\C)$, and completely determines the type of $G_1$.

\subsection{Application of the theory when $G_1$ is not of  dihedral type}

If $G_1$ is finite but not of dihedral type,  then  the main
point in  our  proof of Theorem \ref{abfiniteG1}  will be   to show that
equation $L_4 (v) = 0$ does not have rational solutions. This is why we need to
compute the exponents of $L_4$ at the singularities.

\begin{lemma}
  \label{expL4}
With the notation $2 \varepsilon_i =\{2 a, 2 b\}$, for $ \varepsilon_i=\{a,b\}$,  the respective exponents of $L_4$ at $z\in\{0,1,\infty\}$ are the  following
  \[ \{1, 2, 2 \varepsilon_0 \}, \{1, 2, 2 \varepsilon_1 \}, \{- 1, - 1, 2
     \varepsilon_{\infty} \}. \]
\end{lemma}
\begin{proof}
  If $\varepsilon_i =\{a , b\}$ are the exponents of $L_2$ at the singularity
  $i \in \{0, 1, \infty\}$, then the exponents of $L_3 = L_2^{\circledS 2}$ at the
  same singularity are $\{a + b , 2 a , 2 b\}$. Since at $z = 0$ and $z = 1$,
  $a + b = 1$, and $a + b = - 1$ at $z = \infty$, this gives the exponents of
  $L_3$.

  Let  $\chi_3$ and $\chi_4$ be the characteristic polynomials of the equations $L_3=0$ and $L_4=0$, respectively.

  In a neighbourhood  of $z = 0$ we have the following. If \[ L_3 (z^{\rho}) = \chi_3 (\rho) z^{\rho - 3}
  + \cdots, \] then
  \[ L_4 (z^{\rho}) = (\chi_3 (\rho) z^{\rho - 2} + \cdots)' = (\rho - 2)
     \chi_3 (\rho) z^{\rho - 3} + \cdots . \]
  So $\chi_4 (\rho) = (\rho - 2) \chi_3 (\rho)$.

 In a neighbourhood  of  $z = 1$, we obtain a similar result thanks to the
  formula $z (z - 1) =  (z - 1)^2 + (z - 1)$.

In a neighbourhood    $z = \infty$ we have the following.  If
 \[
  v = x^{\rho} + \cdots =
  {z^{-\rho}} + \cdots,
 \]
  then  the first term of $v'''$ is proportional  to $x^{\rho +
  3}$. So, we have
  \begin{eqnarray*}
    L_4 (x^{\rho}) & = & ( \frac{1}{x^2} \chi_3 (\rho) x^{\rho + 3} +
    \cdots)'\\
    & = & (\chi_3 (\rho) x^{\rho + 1} + \cdots)'\\
    & = & (\chi_3 (\rho) \frac{1}{z^{\rho + 1}} + \cdots)'\\
    & = & - (\rho + 1) (\chi_3 (\rho) x^{\rho + 2} + \cdots)
  \end{eqnarray*}
  Hence, up to the sign, $\chi_4 (\rho) = (\rho + 1) \chi_3 (\rho)$.

  Therefore, at $z = 0$ and $z = 1$, the exponents of $L_4$ are those of $L_3$
  together with $\rho = 2$. At $z = \infty$, the exponents of $L_4$ are those
  of $L_3$ together with $\rho = - 1$.
\end{proof}

\begin{corollary}
  \label{absolL4}For all the rows in Table \ref{tablefiniteG},
  except maybe for rows 5 and 6, the equation $L_4=0$ does not have non-zero rational solutions. In
  particular, when $G_1 =\mathcal{G}(k, \lambda)$ is finite but not dihedral,
  $L_4$ does not have non-zero rational solutions.
\end{corollary}
\begin{proof}
  From Table \ref{tablefiniteG} and Lemma \ref{expL4}, we see that for each
  possible case, the exponents of $L_4$ \ at $z = 0$ and $z = 1$ are greater
  or equal to zero. So, if we look for a rational solution $v$ of $L_4 = 0$,
  $v$ must be a polynomial of degree equal to the opposite of one exponent at
  the infinity. Therefore, $\deg (v) \in \{1, - 2 \varepsilon_{\infty} \}$. Hence,
  $\deg (v)$ must be equal to 1, unless maybe, $- 2 \varepsilon_{\infty}$
  contains an integral number $\geq 2$. But for all the rows of Table
  \ref{tablefiniteG}, the set $2 \varepsilon_{\infty}$ does not  contain any integral
  number,  so the possible polynomial to check are of the form $v = z + d$. We have
  \begin{equation*}
    L_4 (z + d)  =  [z (z - 1) L_2^{\circledS 2} (z+b)]'= -2[z (z - 1) F(z)]',
 \end{equation*}
where
\begin{equation*}
 F(z)=- \frac{1}{2} L_2^{\circledS 2} (z+d)=  2 r + r' (z + d).
\end{equation*}
  Thus, \[
         L_4 (z + d) = 0 \ \Longleftrightarrow\  F (z)
  = \frac{c}{z (z - 1)}
        \]
 for a certain  $c \in\C$. Let us  study
  the behaviour of $F (z)$ around $z = 0$ and $z = 1$. From now, we assume
  that we are not in the cases of rows 5 and 6, in particular $|k| \geq
  3$.

  {{Around $z = 0$}}, $r (z)_0 = \frac{a}{z^2}$ with $a = \frac{(1 /
  k)^2 - 1}{4} \neq 0$ \ Therefore, $r (z)'_0 = \frac{- 2 a}{z^3}$. So, if $d
  \neq 0$, $F (z)_0 = \frac{- 2 ad}{z^3}$. This is incompatible with $F (z) =
  \frac{c}{z (z - 1)}$. Hence, we must check this equation with $d = 0$ (i.e.,
  with $F (z) = 2 r + zr'$).

  {{Around $z = 1$}}, $r (z)_1 = \frac{- 3}{16 (z - 1)^2}$
  therefore, $F (z)_1 = \frac{3}{8 (z - 1)^3}$, and this is still incompatible
  with $F (z) = \frac{c}{z (z - 1)}$.

  Thus, for all the rows except maybe for rows 5 and 6, $L_4 (z + d) \neq  0$,
  and $L_4 = 0$ does not have a non-zero rational solution.
\end{proof}

\begin{proof}[Proof of Theorem \ref{abfiniteG1} for  $G_1$ finite but
  not dihedral.]

Let us   assume  that $G_1$ is finite and is not of dihedral type. This corresponds to
  cases of Table \ref{tablefiniteG}, whose row numbers are  greater than $7$.
  From Corollary \ref{absolL4}, $L_4 = 0$ does not have non-zero  rational solutions.
  Therefore, from Proposition \ref{abeuler},  there exists a non-zero $ x \in V = \mathrm{Sol}
  (L_2 )$,  such that $\int sx^2 \not\in F_1$. So, from Proposition
  \ref{abfinite}, $\forall x \in V\setminus\{0\} , \int sx^2 \not\in F_1$ since $G_1$ is not of
  dihedral type. As a consequence, from Proposition \ref{abcriteredim3},
  $G_3^{\circ}$ is not Abelian and we can conclude thanks to Remark~\ref{jordanredId}.
\end{proof}

\subsection{Application when $G_1$ is of dihedral type}

We have to investigate the cases appearing in row 5, 6 and 7 in Table~\ref{tabledihedral}, which for the convenience  of the reader, we give in Table~\ref{tabledihedral}.
\begin{table}[h]
 \caption{\label{tabledihedral} Cases when  $\mathcal{G}(k , \lambda)$ is of dihedral type.}
  \begin{tabular}{|c|c|c|c|c|}
   \hline
row & k & $\lambda$ & Exponents of $L_2$ at $\{0, \infty\}$ &  $\mathcal{G}(k , \lambda)$ \\
    \hline
    5 & 1 & $0$ & $(0, 1), (- 1 / 4, - 3 / 4)$ & cyclic-dihedral\\
    \hline
    6 & -1 & $1$ & $(0, 1), (- 1 / 4, - 3 / 4)$ & cyclic-dihedral\\
    \hline
    7 & $|k| \geq 3$ & $\frac{1}{2} ( \frac{k - 1}{k} + p (p + 1) k$) &
    $\varepsilon_0$, $(  (2 p - 1)/4, -  (2 p + 3)/4)$ &
    dihedral\\
    \hline
  \end{tabular}
\end{table}

We  follow the  strategy applied above. That is, we prove that $G_3^{\circ}$ is not
Abelian because all the integrals $\int sx^2$ are not algebraic. What is more
difficult here is that we cannot deduce this fact from the existence of one
particular non-algebraic integral. We  begin with the simple cases of rows
5 and 6. Next we consider the case of row 7 which is more technical.

\subsubsection{The case of rows 5 and 6}

From Table \ref{tabledihedral}, the common Riemann scheme of $L_2$ is
\[ P \left\{ \begin{matrix}
     0 & 1 & \infty\\
     0 & 1 / 4 & - 1 / 4\\
     1 & 3 / 4 & - 3 / 4
   \end{matrix}\; z \right\} . \]
A basis of solutions is therefore
\[ x_1 = (z - 1)^{1 / 4}, \quad x_2 = (z - 1)^{3 / 4} . \]
Since $x_1 x_2 = z - 1 \in\C(z)$, here $G_1$ is cyclic and
isomorphic to $\Z/ 4\Z$.

\begin{proposition}
  \label{ablines10-11} In the cases of rows 5 and 6, we have
  \begin{enumerate}
    \item For a non-zero solution $x$ of $L_2 = 0$, the integral $\varphi =
    \int sx^2$ is not algebraic.

    \item The group $G^{\circ} = G^{\circ}_3$ is not Abelian.
  \end{enumerate}
\end{proposition}

\begin{proof}
  Thanks to Proposition \ref{abcriteredim3}, the second point is a consequence
  of the first one.

  Since arbitrary solution of $L_2=0$ can be written as $x = \alpha x_1 + \beta x_2$, the general
  from of $\varphi$ is
  \begin{eqnarray*}
    \varphi & = & \int \frac{1}{z (z - 1)} (\alpha (z - 1)^{1 / 4} + \beta (z
    - 1)^{3 / 4})^2,\\
    \varphi & = & \alpha^2  \int \frac{dz}{z \sqrt{z - 1}} + \beta^2  \int
    \frac{\sqrt{z - 1}}{z} dz + 2 \alpha \beta \mathrm{Log} (z),\\
    \varphi & = & \alpha^2 \varphi_1 + \beta^2 \varphi_2 + 2 \alpha \beta
    \mathrm{Log} (z)
  \end{eqnarray*}
  Since $G_1 \simeq \Z/ 4\Z$, there exists $\sigma \in G_2$,
  such that $\sigma (x_1) = \rmi x_1$ and $\sigma (x_2) = - \rmi x_2$. As $\sigma
  (\varphi) = \int s \sigma (x^2)$, we have
  \begin{equation*}
    \sigma (\varphi)  =  \int s (\rmi \alpha x_1 - \rmi \beta x_2)^2 =  - \alpha^2 \varphi_1 - \beta^2 \varphi_2 + 2 \alpha
    \beta \mathrm{Log} (z) .
  \end{equation*}
  If $\varphi \in F_1$, then $\sigma (\varphi) \in F_1$, and so $4 \alpha \beta
  \mathrm{Log} (z) = \varphi + \sigma (\varphi) \in F_1$ is algebraic. Therefore
  $\alpha \beta = 0$, and $\varphi$ is  proportional either to $\varphi_1$, or
  to $\varphi_2$. But in those two remaining cases, the Taylor expansion of
  the integrand around $z = 0$ shows that each $\varphi_j$ for $j \in \{1 ,
  2\}$ can be written in the form $\varphi_j = \pm \rmi \mathrm{Log} (z) + f_j
  (z)$, where $f_j (z)$ is holomorphic around $z = 0$. Therefore, $\mathcal{M}_0
  (\varphi_j) = \mp 2 \pi + \varphi_j$, and $\varphi_j$ cannot be algebraic
  since it has an infinite number of conjugates by the iteration of
  $\mathcal{M}_0$.
\end{proof}

\subsubsection{The case of row 7}

Here, from Table \ref{tabledihedral}, $k$ and $p$ are relative integers with
$|k| \geq 3$, and the Riemann scheme of $L_2$ is
\[ P_1 \left\{ \begin{matrix}
     0 & 1 & \infty\\
     \dfrac{1}{2} - \dfrac{1}{2 k} & \dfrac{1}{4} & \dfrac{2 p - 1}{4}\\[0.5em]
     \dfrac{1}{2} + \dfrac{1}{2 k} & \dfrac{3}{4} & \dfrac{- 2 p - 3}{4}
   \end{matrix} \; z \right\} . \]
\begin{proposition}
  \label{abline7} In the case of row 7 we have
  \begin{enumerate}
    \item For a non-zero solution $x$ of $L_2 = 0$, the integral $\varphi =
    \int sx^2$ is not algebraic.

    \item The group $G^{\circ} = G^{\circ}_3$ is not Abelian.
  \end{enumerate}
\end{proposition}

As in Proposition \ref{ablines10-11} above, the second point is a
consequence of the first one. But the proof of the first point is going to be
divided into several steps since it is more technical.

Notice that if we change $k$ to
 $k' = - k$,  or $p$ is to $p' = - p - 1$, then the Riemann scheme of
$L_2$ is not changed. Therefore, to prove Proposition \ref{abline7} it is enough
to consider the cases with $k \geq 3$ and $p \geq 0$.

\paragraph{The group $D_{2 N}^{\dag}$.}

The differences of exponents of $L_2$ are $\Delta_0 = 1 / k$, $\Delta_1 = 1 / 2$, and
$ \Delta_{\infty} = p + 1 / 2$. So,  the reduced exponents differences are $1 /
k $, $ 1 / 2 $ and $1 / 2$. Therefore, from \cite{poole} p.128-129, the projective
Galois group of $L_2$, i.e., the image of $G_1$ in $ \mathrm{PSL}(2
,\C)$, is isomorphic to the dihedral group $D_{2 k}$, which is of
order $2 k$. From Lemma \ref{abkovacic}, $G_1$ is necessarily conjugated to a
finite subgroup of $D^{\dag}$ which is not cyclic. That is, $G_1$ is not
a subgroup of  the diagonal group $D_{\mathrm{iag}} = \{ \bigl[\begin{smallmatrix}
  \zeta & 0\\
  0 & 1 / \zeta
\end{smallmatrix}\bigr], \zeta \in\C^{\ast} \}$.
Let $W =
\bigl[\begin{smallmatrix}
  0 & - 1\\
  1 & 0
\end{smallmatrix}\bigr]\in\mathrm{SL}(2 ,\C) $ be the Weyl matrix. We
have the following properties.
\begin{enumerate}
  \item $D^{\dag} = D_{\mathrm{iag}} \cup WD_{\mathrm{iag}}$,

  \item For all $R \in \mathrm{SL}(2 ,\C)$, $WRW^{- 1} = R^{- 1}$,

  \item For all $D \in D_{\mathrm{iag}}$, $(WD)^2 = W^2 = - \id$, and $WD$
  is conjugated to $W$ by an element of $WD_{\mathrm{iag}} $.
\end{enumerate}
By Property 3, we can assume
that $W \in G_1$.  As $W^2=-\id$,
the diagonal subgroup of $G_1$, i.e., $G_1 \cap D_{\mathrm{iag}}$,
contains $- \id$. Since it is a finite cyclic group, it is of even
order $N$,  for a certain $N \in 2\mathbb{N}^{\ast}$. Therefore, as a subgroup of $D^{\dag}$, the group $G_1$ is generated
by $W$ and by
 a matrix $R_{\zeta} = \bigl[\begin{smallmatrix}
  \zeta & 0\\
  0 & 1 / \zeta
\end{smallmatrix}\bigr]$ where, $\zeta$ is a primitive $N$-th  root of unity. This is
the group $D^{\dag}_{2 N}$ of order $2 N$, whose presentation is
\[
 D^{\dag}_{2 N} = < W , R_{\zeta} |\;W^2 = - \id , \;R_{\zeta}^N =
   \id ,\; WR_{\zeta} W^{- 1} = R_{\zeta}^{- 1} > .
\]
The image of $D^{\dag}_{2 N}$ in $ \mathrm{PSL}_2 (\C)$ is
the dihedral group $D_N = D_{2 k}$, in the considered situation.

If $\{x_1, x_2 \}$ is a basis of $V$ in which the representation of $G_1$ is
$D^{\dag}_{2 N}$, then the actions of $W$ and $R_{\zeta}$ on this basis are  given
by the formulae
\[ \left\{ \begin{matrix}
     W (x_1) & = & x_2\\
     W (x_2) & = & - x_1
   \end{matrix} \right. \mtext{and}\left\{ \begin{matrix}
     R_{\zeta} (x_1) & = & \zeta x_1\\[0.5em]
     R_{\zeta} (x_2) & = & \dfrac{1}{\zeta} x_2
   \end{matrix} \right. . \]
Therefore, $W (x_1 x_2) = - x_1 x_2$, and $R_{\zeta} (x_1 x_2) = x_1 x_2$. So,
$(x_1 x_2)^2 \in K =\C(z)$, and $L = K [x_1 x_2]$ is quadratic over
$K$. The group $\mathrm{Gal} (F_1 / L) = < R_{\zeta} >$ is cyclic of order $N = 2 k$.
Since $x_1$ has $N$ distinct conjugates under $\mathrm{Gal} (F_1 / L)$, we have $F_1 = L
[x_1]$. Moreover, $x_1^N \in L =\C(z) [x_1 x_2]$.

\paragraph{Algebraicity of the general integral $\varphi = \int sx^2$. }

\begin{lemma}
  \label{abdihedral} Let $\{x_1 , x_2 \}$ be a basis of $V$ in
  which the representation of $G_1$ is $D^{\dag}_{2 N}$. Then the following statements hold true.
  \begin{enumerate}
    \item If there exists $x_0 = \alpha x_1 + \beta x_2 \in V$ with $\alpha
    \beta \neq 0$ such that $\varphi_0 = \int sx_0^2 \in F_1$, then for all $x
    \in V$ the general integral $\varphi = \int sx^2 \in F_1$.

    \item $\varphi_1 = \int sx_1^2 \in F_1$ iff $\varphi_2 = \int sx_2^2 \in
    F_1$.

    \item $\varphi_1 = \int sx_1^2 \in F_1$ iff there exist $\phi \in
   \C(z) [x_1 x_2]$ such that $\int sx_1^2 = \phi x_1^2$.

    \item If $\{y_1 , y_2 \}$ is a basis of $V$ such that $y_1 y_2$ is at most
    quadratic over $\C(z)$, then up to a permutation of the indices,
    $y_1$ is proportional to $x_1$ and $y_2$ is proportional to $x_2$.
  \end{enumerate}
\end{lemma}

\begin{proof}
  (1)~Since
\[
 R_{\zeta} (x_0) = \zeta \alpha x_1 + \frac{\beta}{\zeta} x_2,\mtext{and}\int sx_0^2 = \alpha^2 \varphi_1 + 2 \alpha \beta \int sx_1 x_2 +
  \beta^2 \varphi_2 \in F_1,
\]
 we deduce that
  \[ \zeta^2 \int sR_{\zeta} (x_0^2) = \zeta^4 \alpha^2 \varphi_1 + \zeta^2 2
     \alpha \beta \int sx_1 x_2 + \beta^2 \varphi_2 \in F_1 . \]
  Since $N = 2 k \geq 6$, we can find two primitive $N$-th roots of unity
  $\zeta$ and $\zeta'$, such that  $\card\{1,\zeta^2 , \zeta'^2\}=3$. Therefore, we obtain an
  identity of the form
  \[ \begin{bmatrix}
       1 & 1 & 1\\
       \zeta^4 & \zeta^2 & 1\\
       \zeta'^4 & \zeta'^2 & 1
     \end{bmatrix} \begin{bmatrix}
       \alpha^2 \varphi_1\\
       2 \alpha \beta \int sx_1 x_2\\
       \beta^2 \varphi_2
     \end{bmatrix}= \begin{bmatrix}
       f_1\\
       f_2\\
       f_3
     \end{bmatrix}\in F_1^3 ,  \]
  where the  $3 \times 3$ Vandermonde matrix on the left hand side is invertible. It implies that $\varphi_1,\int sx_1 x_2,\varphi_2\in F_1$,  because $\alpha \beta \neq
  0$. Therefore, by Proposition \ref{abfinite}, any general integral
  $\varphi = \int sx^2 \in F_1$.

  (2)~From Proposition \ref{abfinite} again,
\[
\left( \int sx_1^2 \in F_1\right)  \ \Longrightarrow
 \ \left( \int sW (x_1^2) = \int sx_2^2 \in F_1\right).
\]
  (3)~If $\int sx_1^2 \in F_1$, then, as   $F_1 = L [x_1]$,  we have
  \[ \int sx_1^2 = \sum_{i = 0}^{N - 1} \phi_i x_1^i , \]
and this equality implies that
\[ sx_1^2 = \sum_{i = 0}^{N - 1}\left(\phi_i' + i \frac{x_1'}{x_1} \phi_i\right)
     x_1^i. \]
  But, $x_1^N \in L$ implies that $x_1' / x_1 \in L$, and the above formula
  gives an expansion of $sx_1^2$ in the $L$-basis $\{1, x_1, \cdots, x_1^{N
  - 1} \}$. Therefore $\phi_i' + i \frac{x_1'}{x_1} \phi_i = 0$ for $i \neq 2$
  and \[ \phi_2 + 2 \frac{x_1'}{x_1} \phi_2 = s,\] that is, $\int sx_1^2 = \phi_2
  x_1^2$ with $\phi_2 \in L=\C(z)[x_1x_2]$.

  (4)~If $y_1 y_2$ is at most quadratic over $\C(z)$, its orbit
  under $G_1$ contains at most two elements. Looking at the orbit under the
  subgroup generated by $R_{\zeta}$, we deduce that $y_1 y_2$ must be fixed by
  the subgroup of the rotations $R_{\lambda}$ where $\lambda$ ranges over the
  $k=N/2$ roots of unity. Now, let us write
  \[ y_1 = ax_1 + bx_2 \text{ and } y_2 = cx_1 + dx_2 \text{ with } ad - bc
     \neq 0. \]
  We get the following two expressions
  \[
 \left\{ \begin{split}
       y_1 y_2 & =  acx_1^2 + (bc + ad) x_1 x_2 + bdx_2^2,\\
       R_{\lambda} (y_1 y_2) & =  \lambda^2 acx_1^2 + (bc + ad) x_1 x_2 +
       \frac{bd}{\lambda^2} x_2^2.
     \end{split} \right.
 \]
  But from the proof of point 1, it follows that the family $\{x_1^2 , x_1 x_2 ,
  x_2^2 \}$ is $\C$-linearly independent. Therefore,  when $\lambda$ is
  a $k$-th roots of unity, the equality  $R_{\lambda} (y_1 y_2) = y_1 y_2$
  implies that
  \[ (\lambda^2 - 1) ac = 0 = (1 - \frac{1}{\lambda^2}) bd. \]
  As  $k \geq 3$, we deduce that $ac = bd = 0$, and, up to a
  permutation, $y_1$ is proportional to $x_1$ and $y_2$ is proportional to
  $x_2$.
\end{proof}

From the above lemma, $\{x_1,x_2\}$ is a distinguished basis of $V$, and we have to compute it in order to study the algebraicity of $\varphi_1$.

\paragraph{Solutions of $L_2 = 0$ and the Jacobi polynomials.}

The Jacobi polynomials $J_n^{(\alpha, \beta)}(t)$  with parameters
$(\alpha, \beta)$, and $n \in \mathbb{N}$ are defined by the following formulae
\[ J_n^{(\alpha, \beta)} (t) = \frac{(t - 1)^{- \alpha} (t + 1)^{- \beta}}{2^n
   n!}  \frac{d^n}{dt^n} \left( (t - 1)^{\alpha + n} (t + 1)^{\beta + n}
   \right) ,
 \]
see  p.~95 in \cite{poole}.
The polynomial $J_n^{(\alpha, \beta)} (t)$ is  of degree $n$, and  belongs to the
Riemann scheme
\[
P_J \left\{ \begin{matrix}
     - 1 & \infty & 1\\
     0 & - n & 0\\
     - \beta & \alpha + \beta + n + 1 & - \alpha
   \end{matrix}\; t \right\},
 \]
thus it is a solution of the following equation
\begin{equation}
 \label{eq:jac}
(1-t^2) \frac{d^2 w}{dt^2} +[(\beta-\alpha) -( \alpha+\beta+2)t] \frac{d w}{dt} +n( \alpha+\beta+n+1)w =0
\end{equation}
If $\alpha$ and $\beta$ are real and greater than $- 1$, then  the pairing
\[ \langle P,Q \rangle = \int_{- 1}^1 (1 - t)^{\alpha} (1 + t)^{\beta} P (t) Q (t) dt, \]
defines a scalar product on $\mathbb{R}[t]$.  It can be shown, see, e.g., page 97  in \cite{poole},  that the
family $\{J_n^{(\alpha, \beta)} (t)\}_{n \in \mathbb{N}}$ is an orthogonal
basis for this scalar product. From this it can be proved, see Ex.~2.39 on  p. 94
in \cite{ramisvol2},  that the roots of $J_n^{(\alpha, \beta)}
(t)$ are simple and contained in the real interval $] - 1 , 1 [$.

In what follows, we use the Jacobi polynomials with parameters
\[ (\alpha, \beta) = (- 1 / k, 1 / k), \]
and we denote them by $J_n (t)$.

Using  the following change of variable
\begin{equation}
  t = \sqrt{1 - z} \ \Longleftrightarrow z = 1 - t^2, \label{eqz-t}
\end{equation}
the solutions of $L_2 = 0$ can be expressed  in the terms of variable $t$. We have the following.
\begin{lemma}
  \label{abjacobi} Let $k$ and $p$ be natural integers with $k
  \geq 3$, and
  \[ \left\{
\begin{split}
       x_1 & =  (1 - t^2)^{1 / 2} (- t^2)^{1 / 4}  \left( \frac{t + 1}{t - 1}
       \right)^{1 / 2 k} J_p (t),\\
       x_2 & =  \rmi (1 - t^2)^{1 / 2} (- t^2)^{1 / 4}  \left( \frac{t + 1}{t -
       1} \right)^{- 1 / 2 k} J_p (- t).
     \end{split} \right.
\]
  Then,
   $\{x_1, x_2 \}$ is a basis of $V$ in which the
    representation of $G_1$ is $D^{\dag}_{2 N}$, and
    $L =\C(z) [x_1 x_2] =\C(z) ( \sqrt{1 -
    z}) =\C(t)$.
\end{lemma}

Before proving the above lemma, we use it to finish the proof of Proposition~\ref{abline7}.

\paragraph{The integral $\int sx_1^2 \neq \phi x_1^2$.}
From Lemmas \ref{abdihedral} and \ref{abjacobi}, $\int sx_1^2 \in F_1$ iff
there exists $\phi \in\C(z) ( \sqrt{1 - z}) =\C(t)$ such
that $\int sx_1^2 = \phi x_1^2$, or equivalently
\begin{equation}
  \frac{d}{dz} \left( \phi x_1^2) = sx_1^2 . \label{abeqL12int} \right.
\end{equation}
But here all the quantities may be expressed in term of $t = \sqrt{1 - z}$.
Applying the chain rule we obtain
\begin{equation}
  \frac{d}{dt} \left( \phi x_1^2) = - 2 tsx_1^2 . \label{abeqL12int2} \right.
\end{equation}
From Lemma \ref{abjacobi},
\[
x_1^2 = (1 - t^2) (- t^2)^{1 / 2}  \left( \frac{t + 1}{t - 1} \right)^{1 /
   k} J_p^2 (t) = \rmi t (1 - t^2) \left( \frac{t + 1}{t - 1} \right)^{1 / k}
   J_p^2 (t).
 \]
Since  \[ s = \frac{1}{2kz (1 - z)} = \frac{1}{2kt^2 (t^2 -
1)},\]  equation~\eqref{abeqL12int2}) reads
\begin{equation}
  \frac{d}{dt} \left( \phi (t) t (1 - t^2) \left( \frac{t + 1}{t - 1}
  \right)^{1 / k} J_p^2 (t) \right)  =  \frac{1}{k} \left( \frac{t + 1}{t - 1}
  \right)^{1 / k} J_p^2 (t) .  \label{abeqL12int3}
\end{equation}
If we set $\psi (t) = \phi (t) t (1 - t^2) J_p^2 (t)$, then $\phi \in
\C(t) $ iff  $ \psi \in\C(t)$. In terms of $\psi(t)$,   equation~\eqref{abeqL12int}
has the form
\begin{eqnarray}
  \frac{d}{dt} \left( \psi (t) \left( \frac{1 + t}{1 - t} \right)^{1 / k}
  \right) & = & \frac{1}{k} \left( \frac{1 + t}{1 - t} \right)^{1 / k} J_p^2 (t),
  \label{abeqL12int4}\\
  \frac{d \psi}{dt} + \frac{1}{k}  \left( \frac{1}{1 + t} + \frac{1}{1 - t}
  \right) \psi & = & \frac{1}{k} J_p^2 (t) . \label{abeqL12dif}
\end{eqnarray}
We can use the above equation to  study  the local behaviour  of the  function $\psi(t)$. A simple analysis shows that
 if $\psi(t)$ is rational, then it has a simple zero at  $t = \pm 1$.
 Moreover, by the Cauchy theorem, $\psi$ has no pole in $] -
1 , 1 [$. Since $0 < 1 / k \leq 1 / 3$, if $\psi$ is rational, then the real
function
\[
 t \mapsto \psi (t) \left( \frac{1 + t}{1 - t} \right)^{1 / k},
\]
vanishes at $t = \pm 1$. Therefore, integrating \eqref{abeqL12int4},  we get
\[ 0 = \int_{- 1}^1 \frac{1}{k} \left( \frac{1 + t}{1 - t} \right)^{1 / k} J_p^2 (t)
   dt, \]
but it is impossible since the integrand is positive on $] - 1 , 1 [$.

Therefore, $\varphi_1 = \int sx_1^2$ does not belong to $F_1$.

\paragraph{Proof of Proposition \ref{abline7}. }

\begin{proof}
  Since $\varphi_1 = \int sx_1^2 \not\in F_1$,  by Lemma \ref{abdihedral},
  $\varphi_2 = \int sx_2^2 \not\in F_1$, and for all non-zero $x \in V$ the
  general integral $\varphi = \int sx^2$ does not belong to $F_1$. Therefore, by Proposition~\ref{abcriteredim3},
  $G_3^{\circ}$ is not Abelian.
\end{proof}

\paragraph{Proof of Lemma \ref{abjacobi}.}
 We can prove the first part of the  lemma directly by making a change of dependent and independent variables  in the  equation $x''=rx$.   Namely, if we put
\[
  y = x(z(t))=  (1 - t^2)^{1 / 2} (- t^2)^{1 / 4}  \left( \frac{t + 1}{t - 1}
       \right)^{1 / 2 k} w (t), \mtext{where} z=z(t)=1-t^2,
\]
then $w(t)$  satisfies  equation~\eqref{eq:jac} with $\beta=-\alpha=1/k$  and $n=p$.
Also, we can prove this part of the lemma applying successive transformations of Riemann schemes, see  Chapter VI in~ \cite{poole}.

 This implies that  the function
  \begin{equation}
     y_1= (1 - t^2)^{1 / 2} (- t^2)^{1 / 4}  \left( \dfrac{t + 1}{t - 1}
    \right)^{1 / 2 k} J_p (t). \label{xjacobi1}
  \end{equation}
is a solution of $L_2=0$ expressed in $t$ variable. Moreover, it can be easily shown that
  \begin{equation}
   y_2= \mathcal{M}_1 (y_1)  = \rmi (1 - t^2)^{1 / 2} (- t^2)^{1 / 4}  \left( \dfrac{t + 1}{t - 1}
    \right)^{- 1 / 2 k} J_p (- t). \label{xjacobi2}
  \end{equation}
  Hence,
  \begin{equation}
    y_1 y_2 = t (1 - t^2) J_p (t) J_p (- t) \in\C[t] .
  \end{equation}
  Since $\C(t) /\C(z)$ is quadratic, $y_1 y_2$ is at most
  quadratic over $\C(z)$. Therefore, from point 4 of Lemma~\ref{abdihedral} we deduce that, up to a permutation, $y_1$ is proportional to
  $x_1$ and $y_2$ is proportional to $x_2$. Therefore, $\{y_1, y_2 \}$ is a
  basis of $V$ in which the representation of $G_1$ is $D_{2 N}^{\dag}$, and  we can
  call it $\{x_1, x_2 \}$.

  Since $\C(z) [y_1 y_2] = L =\C(z) [x_1 x_2] \subset
 \C(t)$,  and $W (x_1 x_2) = - x_1 x_2$, element $x_1 x_2$ is quadratic over
  $\C(z)$. Thus,  we deduce that
  \[ L =\C(z) [x_1 x_2] =\C(t) , \]
and this finishes the proof.
\qed

\paragraph{Conclusion.}

From this study, it follows that the first three points of Theorem~\ref{abclassification} are proved.

\section{Symmetries in Table~\ref{MR} and potentials of degree $k = \pm 2$}

In this section, we notice  an important  symmetry contained in Table~\ref{MR}.  We  use it  to prove Theorem \ref{abclassification}, for the exceptional cases when
$\deg (V) = k = \pm 2$.

\subsection{Symmetries in Table~\ref{MR}}
Let us recall that the reduced  VE \eqref{eqblock}  depends on two rational functions $r, s\in\C(z)$. The function $r$ is defined by the equations \eqref{eq:rl} and \eqref{eqdeltaexp}; the function $s$ is the following
\[
s =
  \dfrac{1}{2 kz (z - 1)} .
\]

In Table \ref{tableMR} there are symmetries between the rows for which $k$ is
changed into $\tilde{k} = - k$. In fact we have the following.

\begin{proposition}
  \label{symmetrie} If the pair $(k , \lambda)$ is changed into
  the pair $( \tilde{k} , \tilde{\lambda}) = (- k , 1 - \lambda)$, then
  \begin{enumerate}
    \item The pair of function $(r , s)$ is changed into $( \tilde{r} ,
    \tilde{s}) = (r , - s)$.

    \item For all $d \geq 1$ the differential Galois groups $G_d$ and $\tilde{G}_d$ of the subsystems of  the VE
    associated to the  Jordan blocks $B(\lambda,d)$ and $B(\tilde\lambda,d)$ are
    isomorphic.
  \end{enumerate}
\end{proposition}

\begin{proof}
 If $k$ is changed into $\tilde{k} = - k$, then from \eqref{eqRS} the Riemann
  schemes $P$ (resp. $\tilde{P}$) of the equations $x'' = rx$ (resp. $x'' =
  \tilde{r} x$) have the same exponents at $z = 0$ and at $z = 1$.
Now, $P=\tilde P$ iff $\tilde\tau=\pm\tau$.
 From \eqref{eqdeltaexp} this happens  iff $ \tilde{\lambda} = 1 -
     \lambda$.
  Therefore, if $( \tilde{k} , \tilde{\lambda}) = (- k , 1 - \lambda)$, then $P
  = \tilde{P}$, and $\tilde{G}_1 = G_1$.  Moreover, from \eqref{eqRS} again, we have  $(\tilde r, \tilde s)=(r,-s)$.

  Now, let us make the following change of variables
  \[ \tilde{x} = - x, \quad\tilde{y} = y,  \quad  \tilde{u} = - u  , \]
  in the system
  \[
       x'' = rx,\quad
       y'' = ry + sx, \quad
       u'' = ru + sy.
\]
  We  can easily obtain
\[
       \tilde{x}'' = \tilde{r}  \tilde{x},\quad
       \tilde{y}'' = \tilde{r}  \tilde{y} + \tilde{s}  \tilde{x}, \quad
       \tilde{u}'' = \tilde{r}  \tilde{u} + \tilde{s}  \tilde{y}
  \]
  By considering the first two equations of both systems we see that the two
  Picard-Vessiot extensions $F_2 /\C(z)$ and $\tilde{F}_2
  /\C(z)$ are equal. So their differential Galois group $G_2$ and
  $\tilde{G}_2$ coincide. Similarly, by considering the three equation of both systems  we have
  \[ F_3 = \tilde{F}_3 \text{ and } G_3 = \tilde{G}_3 . \]
  This arguments are can be obviously generalised  for any Jordan block of size $d
  \geq 3$.
\end{proof}

As a consequence, Table \ref{tableMR} remains stable for the involutive pairing $(k,
\lambda) \leftrightarrow ( \tilde{k}, \tilde{\lambda})$. For example, for
rows 2, 3 and 4, we have \
\[ \lambda (k , p) + \lambda (- k , 1 - p) = 1. \]
So, if $\lambda = \lambda (k , p)$ then $\tilde{\lambda} = \lambda (- k , 1 -
p)$.

\subsection{The case $k = \pm 2$}

\begin{proposition}
  Let $V (q)$ be a homogeneous potential of degree $k = \pm 2$. Then at an arbitrary
  PDP, the connected component
  $G (\mathrm{VE}_t)^{\circ} \simeq G (\mathrm{VE}_z)^{\circ}$ is Abelian.
\end{proposition}

\begin{proof}
  Let us assume that $k = 2$. The VE (\ref{eqVEt}) $\ddot\eta  = - \varphi^{k -
  2} V'' (c) \eta $ reduces to the following linear differential system
  with constant coefficients
  \[ \ddot \eta  = - V'' (c) \eta  . \]
  Let $F /\C(\varphi (t), \dot\varphi (t))$ be the Picard-Vessiot
  extension associated to this system. It is generated over
  $\C(\varphi (t), \dot\varphi(t))$ by the entries of a $n \times n$
  matrix $\Xi (t) = \exp (St)$, where, $S$ is a constant matrix such that
  \begin{equation}
    S^2 = - V'' (c) . \label{eqScarre}
  \end{equation}
  Since it is always possible to extract a square root of a complex matrix,
  \eqref{eqScarre} has a solution whose spectrum consists of numbers $\mu_i$
  with $\mu_i^2 = - \lambda_i$, where  the $\lambda_i$ belong to the spectrum of
  $V'' (c)$. By considering the Jordan decomposition $S = D + N$ of $S$ with
  $D$ conjugated to $\mathrm{diag} (\mu_1, \ldots, \mu_n)$, the entries of $\Xi
  (t)$ are polynomial in $t$ combinations of the exponential $\exp (\mu_i t)$.

  Since the hyperelliptic equation (\ref{eqelliptic}) is now
  \[ \dot\varphi (t)^2 + \varphi (t)^2 = 1 \Rightarrow \ddot\varphi(t) = - \varphi
     (t), \]
  the associated ground field is $\C(\varphi (t), \dot\varphi (t))
  =\C(\exp (\rmi t))$. Therefore,  the connected component $G
  (\mathrm{VE}_t)^{\circ}$ is either a torus,  or the direct product of a
  torus and $G_{\mathrm{a}}$. The latter case  happens only if some of the above mentioned
  polynomials appearing inside $\Xi (t)$ are not constant. In both cases $G
  (\mathrm{VE}_t)^{\circ}$ is Abelian, and the same happens for the connected
  component $G (\mathrm{VE}_z)^{\circ}$, by Proposition \ref{VEtVEz}.

  As a consequence for any system of the form \eqref{eqblock}, corresponding
  to a Jordan block of size $d \geq 1$ with $k = 2$, the connected
  component $G_d^{\circ}$ is Abelian. Moreover, this result is independent of
  the value of the eigenvalue $\lambda$.

  Now, let $\tilde k=-2$. Over the ground field $\C(z)$,
  the VE \eqref{eqVEredz}, can be written as a direct sum of $m$ systems of
  the form \eqref{eqblock}, corresponding to Jordan blocks of sizes $d_i$ and
  eigenvalue $\tilde{\lambda}_i$. If we denote by
  $\tilde{G}_{d_i}$ their respective Galois groups for $1 \leq i
  \leq m$, then  from Section 1.4 we know that  $G (\mathrm{VE}_z)^{\circ}$ is  an
  algebraic subgroup of the direct product
  \[ \tilde{G}_{d_1}^{\circ} \times \cdots \times \tilde{G}^{\circ}_{d_m} . \]
  But from the above principle of symmetries, each $\tilde{G}_{d_i} \simeq
  G_{d_i}$, where $G_{d_i}$ is the Galois group of system (\ref{eqblock})
  corresponding to  Jordan blocks of size $d_i$ and eigenvalue $\lambda_i =
  1 - \tilde{\lambda}_i$, with $k = 2$. Since each $G_{d_i}^{\circ}$ is
  Abelian, so does $G (\mathrm{VE}_z)^{\circ}$ and $G (\mathrm{VE}_t)^{\circ}$.
\end{proof}

\section{About the applications of Theorem \ref{abclassification}}

From now,  $n$ and $k$ are fixed integers with $n \geq 2$ and $k\in\Z^\star$, $c \in\C^n \backslash \{0\}$ is a fixed non-zero complex vector.
In $\C^n$ we define the following pairing
\begin{equation*}
 \langle x, y\rangle:=\sum_{i=1}^nx_iy_i, \mtext{where} x=(x_1, \ldots, x_n)\in\C^n, \quad
 y=(y_1, \ldots, y_n)\in\C^n.
\end{equation*}
Our aim in  Section 6.1 is to show the existence of a great
amount of homogeneous polynomial potentials of degree $k$ such that $c$ is a
PDP of $V$ and $V'' (c) = A$ is a $n \times n$ symmetric matrix as general as
possible. As a consequence, there are  a lot of potentials such that $V'' (c)$
is not diagonalizable. Next, in Section 6.3, we find an  explicit condition for the integrability which does not involve the eigenvalues of the Hessian.

\subsection{From polynomial potential to symmetric matrices}

Here we assume that $k\geq3$, and we  consider the following sets
\begin{align*}
  &R_{n, k}  = \defset{V (q) \in\C[q_1, \ldots, q_n]}{
      V \mtext{is homogeneous, and} \deg (V) = k},\\
&R_{n, k} (c)  =  \defset{V (q) \in R_{n, k}}{V' (c) = c},\\
 &\mathrm{Sym}_n  =  \defset{A \in M_n (\C)}{ A^T = A},\\
&\mathrm{Sym}_{n, k} (c)  =  \defset{A \in \mathrm{Sym}_n}{Ac = (k - 1) c}.
\end{align*}
All these sets are affine spaces of respective complex dimensions
\begin{align*}
  &\dim R_{n, k} =\dbinom{n + k - 1}{ n - 1}, &\dim R_{n, k} (c) =\dbinom{n + k - 1}{ n - 1}-n, \\
 &\dim \mathrm{Sym}_n=\dbinom{n+1}{2}, &\dim \mathrm{Sym}_{n, k} (c) =\dbinom{n+1}{2}-n.
\end{align*}
Now, the Hessian map
\[ h : R_{n, k} (c) \rightarrow \mathrm{Sym}_n, \quad V \mapsto V'' (c) \]
is an affine morphism whose image is contained in $\mathrm{Sym}_{n, k} (c)$.
Indeed, from the Euler identity,
\[ V' (c) = c \ \Longrightarrow\ V'' (c) c = (k - 1) c. \]
More precisely we have the following property whose proof follows from computations of dimensions.
\begin{proposition}
  The image of the Hessian map coincides with $\mathrm{Sym}_{n, k} (c)$. In
  other words, if $n \geq 2$ and $k \geq 3$, then for an arbitrary complex
  symmetric matrix satisfying $Ac = (k - 1) c$, there exists a homogeneous
  polynomial potential $V$ of $n$ variables and degree $k$ such that $c$ is
  PDP of $V$ and $V'' (c) = A$.
\end{proposition}

\subsection{Non diagonalizable complex symmetric matrices}

Let us assume that  $k\in\Z^*$, we show that there are a lot of non-diagonalizable symmetric
matrices belonging  to  the space
$\mathrm{Sym}_{n, k} (c)$.   The most reachable ones belong to
\[ \mathrm{Spec}_{k - 1} :=\defset{A \in \mathrm{Sym}_{n, k} (c)}{\mathrm{Spec} (A) =\{k
   - 1\}}, \]
which is the set of matrices such that $\lambda = k - 1$ is the only
eigenvalue of $A$. Indeed, any such $A$ is either equal to $(k - 1) \cdot
\id$, or non-diagonalizable.
\begin{proposition}
\label{pro:62}
  With the notations above we have:
  \begin{enumerate}
    \item If $n \geq 3$, then  the space $\mathrm{Sym}_{n, k} (c)$ contains non-diagonalizable matrices.

    \item For $n = 2$, the space  $\mathrm{Sym}_{2, k} (c)$ contains  non
    diagonalizable matrices iff $c$ is isotropic, i.e.,  $\langle c, c \rangle = 0$.
\item Moreover, when $n=2$ and $c$ is isotropic,  $\mathrm{Sym}_{2, k} (c)= \mathrm{Spec}_{k - 1}$.
  \end{enumerate}
\end{proposition}

\begin{proof}
  By triangularizing a $n \times n$ matrix $A$, we see that $\lambda=k - 1$ is
  its only eigenvalue iff $A$ satisfies the following $n$ algebraic equations
  \begin{equation}
   \frac{1}{n} \tr A^p = (k - 1)^p \mtext{ for } 1 \leq p \leq
    n. \label{traceA}
  \end{equation}
  If $A \in \mathrm{Sym}_{n, k} (c)$, then $\lambda = k - 1$ is one of the eigenvalues of
  $A$, so we only need the $(n - 1)$ first equations of \eqref{traceA} to
  ensure that its $(n - 1)$ remaining eigenvalues coincide with $k - 1$. This
  proves that $\mathrm{Spec}_{k - 1}$ is an algebraic affine subset of
  $\mathrm{Sym}_{n, k} (c)$, whose dimension satisfies
  \[ \dim \mathrm{Spec}_{k - 1} \geq \dim \mathrm{Sym}_{n, k} (c) - (n - 1) =\binom{n+1}{2}
     - 2 n + 1. \]
  Hence $\dim \mathrm{Spec}_{k - 1} \geq 1$, as soon as $n \geq 3$.
  This proves point 1.

  Let $n = 2$, and  $A \in \mathrm{Sym}_{2, k} (c)$. The line
  $(\C c)^{\bot}$ of vectors of $\C^2$ which are orthogonal
  to $c$ is globally left invariant by $A$. Therefore, if
  $(\C c)^{\bot} \neq \text{$\C c$}$, then $A$ is
  diagonalizable in a basis $(c, v)$, where $v \in (\C c)^{\bot}$. So,
  if $\mathrm{Sym}_{2, k} (c)$ contains a matrix which is not diagonalizable, then we have $(\C c)^{\bot} = \text{$\C c$}$. That
  is, $c$ is  isotropic. Conversely, let us assume that $c$ is isotropic.
Let $A\in\mathrm{Sym}_{2, k} (c)$. If $\lambda\neq k-1$ is another eigenvalue of $A$, then $\ker(A-\lambda\id)$ is a line which is orthogonal to $\C c$ and different from it. This is impossible since  $(\C c)^\perp=\C c$. Therefore $k-1$ is the only possible eigenvalue of each matrix belonging to  $\mathrm{Sym}_{2, k} (c)$. Hence, $\mathrm{Sym}_{2, k} (c)= \mathrm{Spec}_{k - 1}$, and point 3 is proved.  Since
\[
 \dim \mathrm{Sym}_{2, k} (c)= 1 =\dim\mathrm{Spec}_{k - 1},
\]
except for the matrix $(k-1)\cdot \id$, any other matrix of $\mathrm{Sym}_{2, k} (c)$ is not diagonalizable. This proves point 2.
\end{proof}

\subsection{New necessary condition for integrability}
Here we focus our attention on some potentials admitting isotropic PDP.   Altought the eigenvalues of the  Hessian $V''(c)$ does not give  any obstacle to the integrability,  we exhibit a new one.
  In the following  we set $c_0=(1,\rmi)$.
\begin{proposition}
 Let $V(q)=V(q_1,q_2)$  be a two degrees of freedom homogeneous potential of the following form
\begin{equation*}
 V(q):=(q_1^2+q_2^2)W(q)
\end{equation*}
where $W(q)$ is a homogeneous function  with
\begin{equation*}
 \deg W\in\Z\setminus\{ -4,-2,-1,0\}, \quad W(c_0)\in\mathbb{P}^1\setminus\{0,\infty\}.
\end{equation*}
If the Hamiltonian system associated with this potential is completely integrable, then
\begin{equation*}
\rmi\pder{W}{q_1}(c_0) +\pder{W}{q_2}(c_0)= 0.
\end{equation*}
\end{proposition}
\begin{proof}
 If $k=\deg V=2 +\deg W$, then $k \in \Z\setminus\{-2,0,1,2\}$, i.e., we have either $|k|\geq 3$, or $k=-1$. From
$V(q)=(q_1^2+q_2^2)W(q)$ we get
\begin{equation*}
 V'(c_0) =2 W'(c_0)\cdot c_0.
\end{equation*}
So, $V'(c)=c$ for $c= \mu c_0$, where $2\mu^{k-2}W(c_0)=1$. Hence, according to point 3 of Proposition~\ref{pro:62}, it follows that $V''(c)\in \mathrm{Spec}_{k - 1}$, i.e., $\lambda=k-1$ is the only eigenvalue of $V''(c)$. Hence, $G_1^\circ \simeq G_\mathrm{a}$, and the potential satisfies  the conditions  appearing in the row 2 or 3 of Table~\ref{tableMR}. Thus Theorem~\ref{MR} does not give any obstacles for the integrability of $V$.  Now, Theorem~\ref{abclassification} gives an obstacle iff $V''(c)$ is not diagonalizable. This happens iff $V''(c)\neq (k-1)\cdot \id$, i.e., iff
\begin{equation*}
 \frac{\partial^2 V}{\partial q_1 \partial q_2}(c)\neq 0 \quad \Longleftrightarrow \quad \frac{\partial^2 V}{\partial q_1 \partial q_2}(c_0)\neq 0
\end{equation*}
But a direct computation shows that the last condition is equivalent  to the following one
\begin{equation*}
 \rmi\pder{W}{q_1}(c_0) +\pder{W}{q_2}(c_0)\neq 0.
\end{equation*}
\end{proof}
Let $V(q)= (q_1^2+q_2^2)W(q)$, with $ \deg W\in\Z\setminus\{ -4,-2,-1,0\}$, and  $W(c_0)\in\mathbb{P}^1\setminus\{0,\infty\}$. The condition
\begin{equation*}
 \rmi\pder{W}{q_1}(c_0) +\pder{W}{q_2}(c_0)= 0,
\end{equation*}
is therefore a non-trivial condition for the integrability of $V$. For example, if $W(q)=aq_1+bq_2$, then $V$ is integrable if $W(q) = \alpha( q_1 -\rmi q_2)$, $\alpha\in\C$. In this case $ V $ is indeed integrable with the additional first integral
\begin{equation}
 F=\rmi p_1^2 +6p_1p_2 -5\rmi p_2^2 +8\alpha q_2(q_1-\rmi q_2)^2.
\end{equation}
\section*{Acknowledgements}
The first author wishes to thank University of Zielona G\'ora for the excellent conditions during his short visit in February 2008.

For  the second author this research  was supported by grant No. N N202 2126 33 of Ministry of Science and Higher Education of Poland.

\end{document}